\input lanlmac
\def\href#1#2{{#2}}
\def\hhref#1{{#1}}
\input epsf.tex

\overfullrule=0mm

\newcount\figno
\figno=0
\def\fig#1#2#3{
\par\begingroup\parindent=0pt\leftskip=1cm\rightskip=1cm\parindent=0pt
\baselineskip=11pt
\global\advance\figno by 1
\midinsert
\epsfxsize=#3
\centerline{\epsfbox{#2}}
\vskip 12pt
{\bf Fig.\ \the\figno:} #1\par
\endinsert\endgroup\par
}
\def\figlabel#1{\xdef#1{\the\figno}}
\def\encadremath#1{\vbox{\hrule\hbox{\vrule\kern8pt\vbox{\kern8pt
\hbox{$\displaystyle #1$}\kern8pt}
\kern8pt\vrule}\hrule}}


\def\IR{\relax{\rm I\kern-.18em R}}
\font\cmss=cmss10 \font\cmsss=cmss10 at 7pt

\font\cmss=cmss10 \font\cmsss=cmss10 at 7pt
\def\IZ{\relax\ifmmode\mathchoice
{\hbox{\cmss Z\kern-.4em Z}}{\hbox{\cmss Z\kern-.4em Z}}
{\lower.9pt\hbox{\cmsss Z\kern-.4em Z}}
{\lower1.2pt\hbox{\cmsss Z\kern-.4em Z}}\else{\cmss Z\kern-.4em Z}\fi}
\def\IN{\relax{\rm I\kern-.18em N}}
\def\b{\circ}
\def\n{\bullet}

\def\gbbbb{\Gamma_4^{\hbox{$\scriptstyle \b \b$}\kern -8.2pt
\raise -4pt \hbox{$\scriptstyle \b \b$}}}
\def\gnnnn{\Gamma_4^{\hbox{$\scriptstyle \n \n$}\kern -8.2pt  
\raise -4pt \hbox{$\scriptstyle \n \n$}}}
\def\gnnnnnn{\Gamma_6^{\hbox{$\scriptstyle \n \n \n$}\kern -12.3pt
\raise -4pt \hbox{$\scriptstyle \n \n \n$}}}
\def\gbbbbbb{\Gamma_6^{\hbox{$\scriptstyle \b \b \b$}\kern -12.3pt
\raise -4pt \hbox{$\scriptstyle \b \b \b$}}}
\def\gbbbbc{\Gamma_{4\, c}^{\hbox{$\scriptstyle \b \b$}\kern -8.2pt
\raise -4pt \hbox{$\scriptstyle \b \b$}}}
\def\gnnnnc{\Gamma_{4\, c}^{\hbox{$\scriptstyle \n \n$}\kern -8.2pt
\raise -4pt \hbox{$\scriptstyle \n \n$}}}
\def\Rbud#1{{\cal R}_{#1}^{-\kern-1.5pt\blacktriangleright}}
\def\Rleaf#1{{\cal R}_{#1}^{-\kern-1.5pt\vartriangleright}}
\def\Rbudb#1{{\cal R}_{#1}^{\circ\kern-1.5pt-\kern-1.5pt\blacktriangleright}}
\def\Rleafb#1{{\cal R}_{#1}^{\circ\kern-1.5pt-\kern-1.5pt\vartriangleright}}
\def\Rbudn#1{{\cal R}_{#1}^{\bullet\kern-1.5pt-\kern-1.5pt\blacktriangleright}}
\def\Rleafn#1{{\cal R}_{#1}^{\bullet\kern-1.5pt-\kern-1.5pt\vartriangleright}}
\def\Wleaf#1{{\cal W}_{#1}^{-\kern-1.5pt\vartriangleright}}
\def\Cleaf{{\cal C}^{-\kern-1.5pt\vartriangleright}}
\def\Cbud{{\cal C}^{-\kern-1.5pt\blacktriangleright}}
\def\Crleaf{{\cal C}^{-\kern-1.5pt\circledR}}

\def\Encadre#1{\vcenter{\hrule\hbox{\vrule\kern8pt\vbox{\kern8pt#1\kern8pt}
\kern8pt\vrule}\hrule}}
\def\encadremath#1{\Encadre{\hbox{$\displaystyle#1$}}}


\magnification=\magstep1
\baselineskip=12pt
\hsize=6.3truein
\vsize=8.7truein
 at 8truept
 at 8truept
 at 10truept

\font\bigrm=cmr12 at 14pt
\centerline{\bigrm Integrability of graph combinatorics}
\centerline{\bigrm via random walks and heaps of dimers}

\bigskip\bigskip

\centerline{P. Di Francesco and E. Guitter}
  \smallskip
  \centerline{Service de Physique Th\'eorique, CEA/DSM/SPhT}
  \centerline{Unit\'e de recherche associ\'ee au CNRS}
  \centerline{CEA/Saclay}
  \centerline{91191 Gif sur Yvette Cedex, France}
\centerline{\tt philippe@spht.saclay.cea.fr}
\centerline{\tt guitter@spht.saclay.cea.fr}

  \bigskip


     \bigskip\bigskip

     \centerline{\bf Abstract}
     \smallskip
     {\narrower\noindent
We investigate the integrability of the discrete non-linear equation governing the 
dependence on geodesic distance of planar graphs with inner vertices of even valences.
This equation follows from a bijection between graphs and blossom trees and is expressed
in terms of generating functions for random walks. We construct explicitly an infinite 
set of conserved quantities for this equation, also involving suitable combinations of 
random walk generating functions. The proof of their conservation, i.e. their eventual independence on 
the geodesic distance, relies on the connection between random walks and heaps of dimers. 
The values of the conserved quantities are identified with generating functions for 
graphs with fixed numbers of external legs. Alternative equivalent choices for the set 
of conserved quantities are also discussed and some applications are presented.\par}

     \bigskip

\Date{06/05}

\nref\BIPZ{E. Br\'ezin, C. Itzykson, G. Parisi and J.-B. Zuber, {\it Planar
Diagrams}, Comm. Math. Phys. {\bf 59} (1978) 35-51.}
\nref\FIK{A. Fokas, A. Its and A. Kitaev, {\it Discrete Painlev\'e equations
and their appearance in quantum gravity}, Comm. Math. Phys. {\bf 142} (1991)
313-344.}
\nref\QGRA{V. Kazakov, {\it Bilocal regularization of models of random
surfaces}, Phys. Lett. {\bf B150} (1985) 282-284; F. David, {\it Planar
diagrams, two-dimensional lattice gravity and surface models},
Nucl. Phys. {\bf B257} (1985) 45-58; J. Ambjorn, B. Durhuus and J. Fr\"ohlich,
{\it Diseases of triangulated random surface models and possible cures},
Nucl. Phys. {\bf B257}(1985) 433-449; V. Kazakov, I. Kostov and A. Migdal
{\it Critical properties of randomly triangulated planar random surfaces},
Phys. Lett. {\bf B157} (1985) 295-300.}
\nref\MAT{See for instance: P. Di Francesco, P. Ginsparg
and J. Zinn--Justin, {\it 2D Gravity and Random Matrices},
Physics Reports {\bf 254} (1995) 1-131, and references therein; see also
B. Eynard, {\it Random Matrices}, Saclay Lecture Notes (2000),
available at {\sl http://www-spht.cea.fr/lectures\_notes.shtml} }
\nref\SCHth{G. Schaeffer, {\it Conjugaison d'arbres
et cartes combinatoires al\'eatoires}, PhD Thesis, Universit\'e 
Bordeaux I (1998); available at \hhref{http://www.lix.polytechnique.fr/\~{}schaeffe/Biblio}.}
\nref\SCH{G. Schaeffer, {\it Bijective census and random
generation of Eulerian planar maps}, Electronic
Journal of Combinatorics, vol. {\bf 4(1)} (1997) R20.}
\nref\CENSUS{J. Bouttier, P. Di Francesco and E. Guitter, {\it Census of planar
maps: from the one-matrix model solution to a combinatorial proof},
Nucl. Phys. {\bf B645}[PM] (2002) 477-499, arXiv:cond-mat/0207682.}
\nref\GEOD{J. Bouttier, P. Di Francesco and E. Guitter, {\it Geodesic
distance in planar graphs}, Nucl. Phys. {\bf B663}[FS] (2003) 535-567, 
arXiv:cond-mat/0303272.}
\nref\VIENNOT{X. Viennot, {\it Heaps of pieces 1: basic definitions and combinatorial lemmas}, 
in G. Labelle and P. Leroux Eds. {\it Combinatoire \'enum\'erative}, Lect. Notes in Math. 
vol {\bf 1234} (1986) 321-350.} 
\nref\BMR{M. Bousquet-M\'elou and A. Rechnitzer, {\it Lattice animals and heaps
of dimers}, Discrete Math. {\bf 258} (2002) 235-274.}
\nref\DGL{P. Di Francesco and E. Guitter, {\it Critical and multicritical
semi-random $(1+d)$-dimensional lattices and hard objects in $d$
dimensions}, J. Phys. A Math. Gen, {\bf 35} (2002) 897-927.}
\nref\MOB{J. Bouttier, P. Di Francesco and E. Guitter, {\it Planar Maps as Labeled Mobiles},
Electronic Journal of Combinatorics, vol. {\bf 11(1)} (2004) R69.}
\nref\CORA{R. Cori and B. Vauquelin, {\it Planar maps are well labeled trees},
Canad. J. Math. {\bf 33 (5)} (1981) 1023-1042; D. Arqu\`es, {\it Les hypercartes 
planaires sont des arbres tr\`es bien \'etiquet\'es}, Discr. Math. {\bf 58}(1) (1986) 11-24;
M. Marcus and G. Schaeffer, {\it Une bijection simple pour les
cartes orientables}, \hhref{http://www.lix.polytechnique.fr/\~{}schaeffe/Biblio/}.}
\nref\ONEWALL{J. Bouttier, P. Di Francesco and E. Guitter, {\it Statistics
of planar maps viewed from a vertex: a study via labeled trees},
Nucl. Phys. {\bf B675}[FS] (2003) 631-660, arXiv:cond-mat/0307606.}
\nref\CONST{M. Bousquet-M\'elou and G. Schaeffer,
{\it Enumeration of planar constellations}, Adv. in Applied Math.,
{\bf 24} (2000) 337-368.}
\nref\RAMA{P. Di Francesco, {\it Geodesic distance in planar graphs: An integrable approach}, 
preprint (2003), to appear in the Ramanujan Journal}
\nref\JM{M. Jimbo and T. Miwa, {\it Solitons and infinite dimensional Lie
algebras}, Publ. RIMS, Kyoto Univ. {\bf 19} No. 3 (1983) 943-1001,
eq.(2.12).}
\nref\GNR{B. Grammaticos, F. Nijhoff and A. Ramani, {\it Discrete Painlev\'e equations},
in {\it The Painlev\'e property, one century later}, R. Conte Ed., CRM series in Math. Phys.
(1999) 413-516.}
\nref\EMB{M. Bousquet-M\'elou, {\it Limit laws for embedded trees. Applications to the 
integrated superBrownian excursion} preprint arXiv:math.CO/0501266.}


\newsec{Introduction}

Graph combinatorics seems to provide a remarkable source for (discrete) integrable systems.
So far, the underlying integrability seemed to be intimately related to the existence
of a matrix model formulation for the counting of these graphs [\xref\BIPZ,\xref\FIK]. 
In the context of matrix model solutions to 2D Quantum Gravity \QGRA, the generating functions for
possibly decorated graphs of arbitrary genus may indeed be interpreted as tau-functions 
of various integrable hierarchies \MAT.

More recently, an alternative purely combinatorial and bijective approach to 
the enumeration of {\it planar} graphs was developed, based on the transformation
of graphs into decorated trees [\xref\SCHth-\xref\CENSUS]. In this context, another remarkable, apparently
unrelated integrable structure was observed, now involving the 
{\it geodesic distance} on the graphs \GEOD.

More precisely, in Ref.\GEOD, special attention was paid to the case of planar graphs 
with vertices of even valence only, and to their generating function $R_n$ with two
distinguished points at geodesic distance at most $n$, with $n$ a non-negative integer. 
The latter  was shown to obey 
a self-consistent {\it master equation} turning into a recursion relation on $n$ for graphs
with bounded valences.
This non-linear recursion relation turned out to be exactly solvable in terms
of discrete soliton-like tau-functions. The precise solution involved
integration constants which may be rephrased into {\it conserved quantities},
thus displaying explicitly the integrability of the master equation.
The knowledge of these conserved quantities provides a powerful tool for
generating explicit expressions for various generating functions of interest
in the graph combinatorial framework. In particular, exploiting the conserved
quantities allows to bypass the quite tedious use of the exact solutions of
Ref.\GEOD.

The aim of this paper is to construct {\it ab initio} a set of conserved
quantities for the master equation in a purely combinatorial setting, and 
to interpret them in the language of graphs. As will be recalled below, the
master equation has a compact expression in terms of discrete random walks.
Our construction relies therefore on properties of partition functions for such
random walks, as well as for ``hard dimers" on a line, both with suitable weights 
involving the $R_n$'s.  Both objects are known [\xref\VIENNOT,\xref\DGL] to be related 
via some boson/fermion 
type inversion relation. The origin of this relation is best understood by relating random
walks to so-called ``heaps of dimers", a boson-like counterpart of the 
mutually excluding, thus fermion-like, hard dimers.

\medskip
The paper is organized as follows. In Section 2, we first recall the existing bijection
between on the one hand planar graphs with two external legs and inner vertices of 
even valence, and on the other hand planted trees with two kinds of leaves, the
so-called blossom trees (Sect.2.1). This bijection is used to re-derive the
master equation for $R_n$ in Sect.2.2. Other possible applications of the
master equation are discussed in Sect.2.3, as well as its integrable nature in
Sect.2.4. Section 3 is devoted to the study of generating functions for
random walks, heaps of dimers, and one-dimensional hard dimer configurations and
to their relationships. In Section 3.1, we first present a useful fundamental relation for
the generating function of random walks, which is then identified with that of a particular
class of heaps of dimers called pyramids. Details of this identification are gathered in
Appendix A. This identification allows to derive a crucial inversion relation involving the
partition function for hard dimers on a segment. In Section 4, we first give explicit formulas 
for the conserved quantities $\Gamma_{2i}(n)$, $i=1,2,\ldots$ in terms of partition functions 
for random walks (Sect.4.1). The actual conservation of these quantities, i.e. their independence of $n$  
whenever the master equation is satisfied, is proved in Sect.4.2 by extensive use of the 
fundamental and inversion formulas of Section 3. Section 5 concerns the graph interpretation 
of the conserved quantities $\Gamma_{2i}(n)$ in terms of $2i$-point functions, i.e.
generating functions for graphs with $2i$ external legs (Sect.5.1). In the 
case of graphs with bounded valences, say up to $2m$, the infinite set of conserved quantities is
reducible to the finite set of the $m-1$ first ones. The corresponding (linear) relations
between their constant values are derived in Sect.5.2. Finally, Sect.5.3 presents
an alternative bijective proof for the conservation of $\Gamma_{2}(n)$, interpreted
as the generating function for graphs with the two legs in the same face. Possible
extensions of this bijective proof to higher order conserved quantities are also discussed,
in view of particularly suggestive identities equivalent to the conservation of $\Gamma_{2i}$
for $i=2,3,4$. In Section 6, we first present two equivalent alternative sets of
conserved quantities, one respecting the ``time-reversal" ($n\to -n$) symmetry of the 
master equation (Sect.6.1) and the other compacted, i.e. involving the least possible
number of $R_n$'s (Sect.6.2). Proofs of their conservation rely on generalized inversion
formulas, all displayed in Appendix B. To illustrate the interest of using conserved quantities,
we give in Sect.6.3 an example of application to the statistics of neighbors of the
external face in pure tetravalent and pure hexavalent graphs. We gather a few
concluding remarks in Section 7.

\newsec{Generating functions for planar graphs: master equation}
\subsec{From planar graphs to blossom trees}
\fig{A two leg-diagram with distinguished incoming and outcoming legs (a). The
incoming leg is adjacent to the external face. All inner vertices have even
valence. Here the geodesic distance between the legs is $1$. This diagram is transformed into
a blossom tree (e) by cutting edges into bud-leaf pairs as explained in the text. 
Buds (resp. leaves) are represented by black (resp. white) arrows. Here 
the cutting procedure requires performing two turns around the graph: (b)$\to$(c)
and (c)$\to$(d). The outcoming leg serves as root for the blossom tree while
the incoming one is replaced by a leaf. In the blossom tree, each inner vertex
of valence $2k$ is adjacent to exactly $k-1$ buds.}{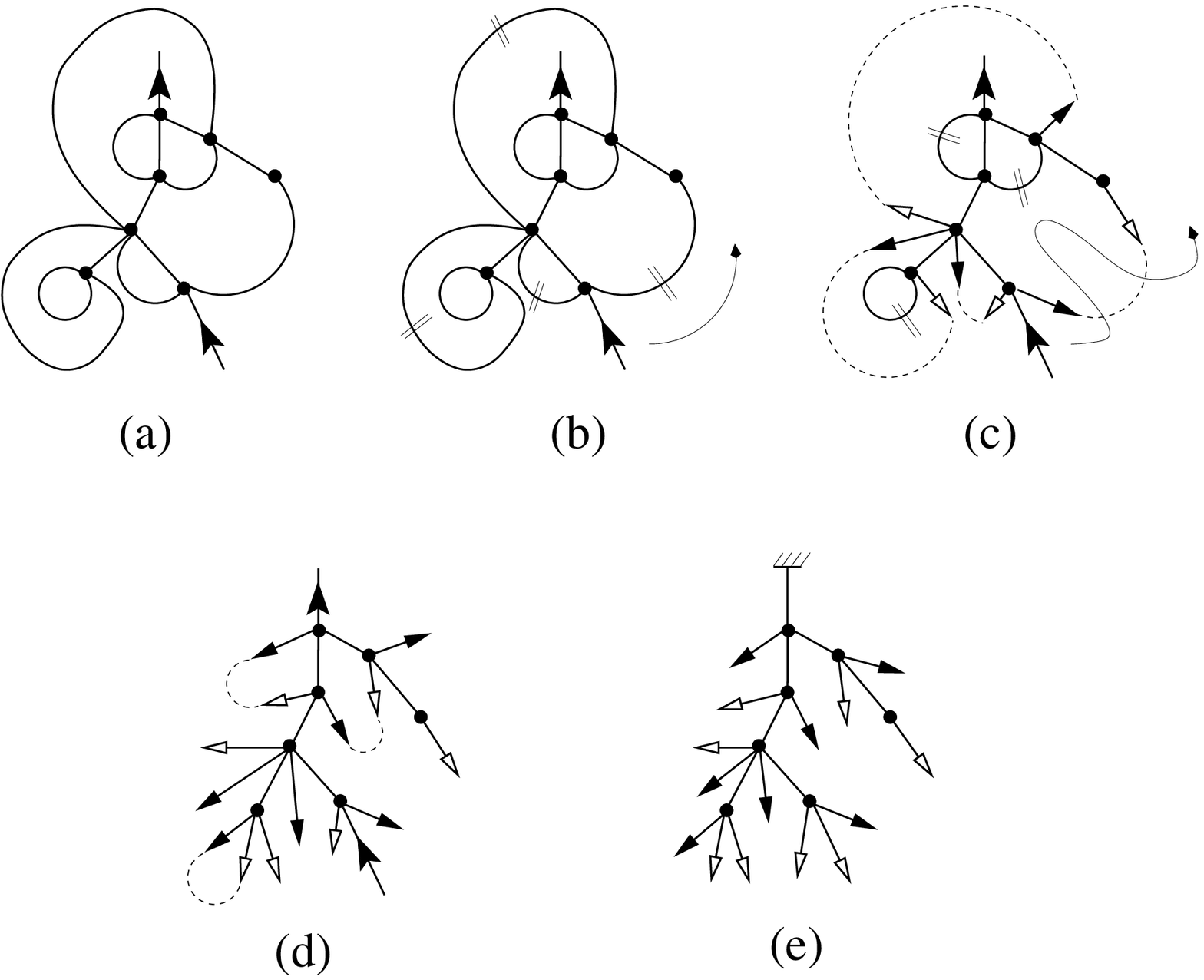}{12.cm}
\figlabel\cutting
Let us first recall the bijection between planar graphs and blossom trees.
More precisely, by planar graphs, we mean here graphs with planar topology,
with inner vertices of {\it even valences} only and with two univalent
endpoints (see Fig.\cutting-(a) for an example). We shall refer to the endpoints
as {\it legs} and to the graphs as {\it two-leg diagrams}. The two 
legs are distinguished into an incoming and an outcoming leg 
and need not be adjacent to the 
same face. We call {\it geodesic distance} between the two legs the minimal number of 
edges to be crossed in a path joining both legs in the plane. In the planar representation, 
we choose to have the incoming leg adjacent to the external face. 

On the other hand, by blossom trees, we mean planted plane trees satisfying
\item{(B1)} all the inner vertices of the tree have even valences;
\item{(B2)} the endpoints of the tree are of two types: buds and leaves;
\item{(B3)} each $2k$-valent inner vertex is adjacent to exactly $k-1$ buds;
\par \noindent (see Fig.\cutting-(e) for an example). In this characterization, 
it is assumed that the tree has at least 
one inner vertex but, by convention, the tree reduced to a single leaf attached to the root
is also considered as a blossom tree. Note that from property (B3), we deduce that 
the total number of leaves in a blossom tree is equal to the total number of buds 
plus one. 

The bijection may be obtained by cutting each two-leg diagram into a blossom tree as follows. 
We first visit successively each edge bordering the external face in counterclockwise direction 
around the diagram (see Fig.\cutting-(b)), starting from the incoming leg. Any
visited edge is cut into a bud-leaf pair if and only if this operation does not disconnect 
the diagram. After one turn, the net result has been to merge
a number of internal faces with the external one. We repeat the process on this 
newly obtained diagram (see Fig.\cutting-(c)) until all faces are merged, resulting into a tree
(see Fig.\cutting-(d)). We finally replace the incoming leg by a leaf and plant the tree at the 
outcoming one (see Fig.\cutting-(e)). As shown in Ref.\SCHth, the resulting tree is a blossom tree
and the above construction establishes a bijection between, on the one hand two-leg diagrams 
and, on the other hand blossom trees.
Note that the inner vertices of the obtained tree are clearly in one-to-one correspondence with 
those of the corresponding diagram and that their valence is preserved. This allows to keep track 
of vertex valences in both pictures.
\fig{The closing procedure from a blossom tree (a) to a two-leg diagram (d).
Each pair of a bud immediately followed counterclockwise by a leaf is glued into
an edge (b). This is iterated (c) until all leaves but one are glued. This remaining
leaf is replaced by an incoming leg while the root becomes an outcoming one. 
The closing procedure here encircles the root once, hence the geodesic
distance between the legs is $1$.}{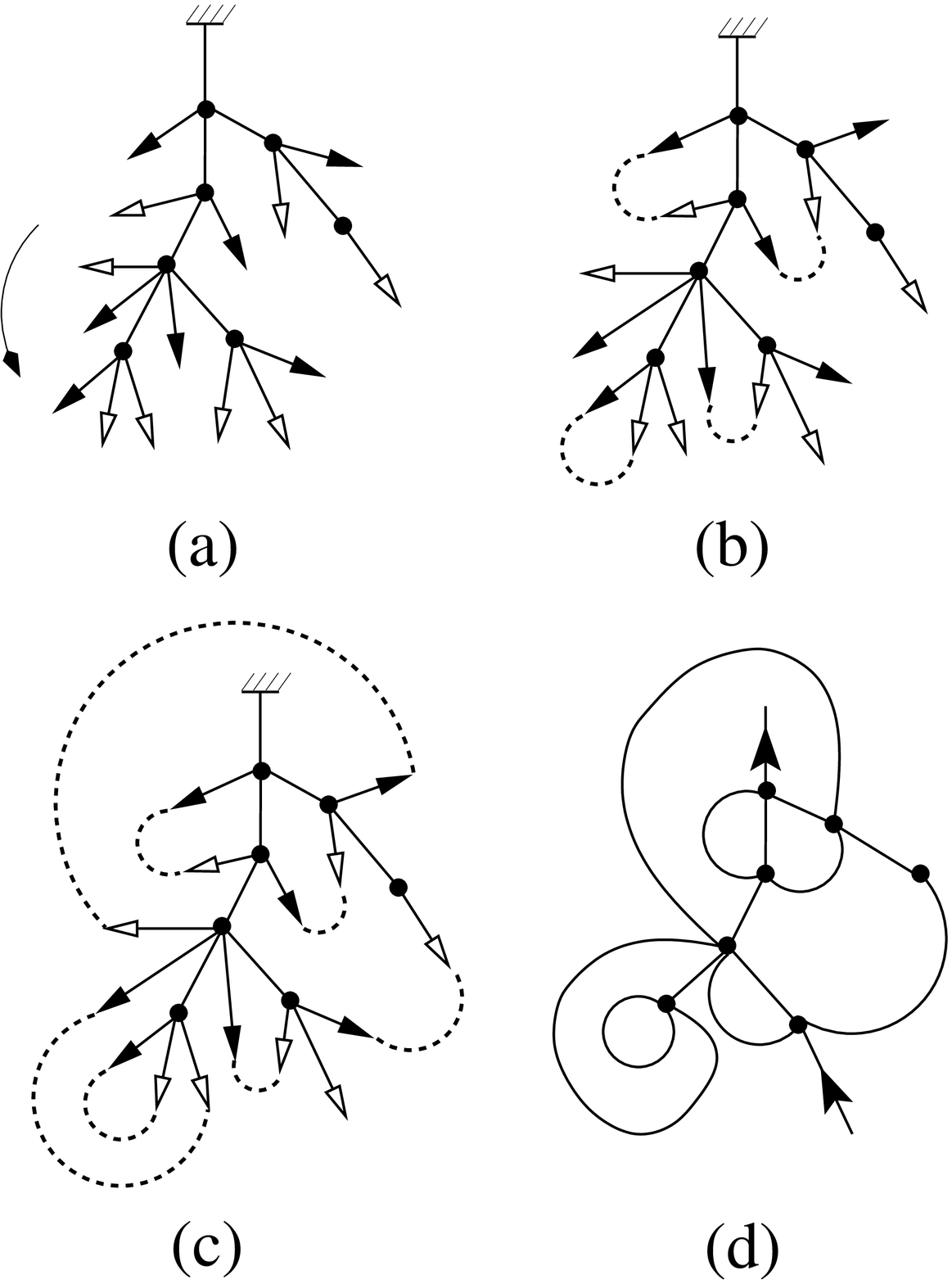}{7.cm}
\figlabel\closing
\fig{The contour and surrounding walks of the blossom tree of Fig.\cutting(e), or
equivalently of Fig.\closing(a). The contour walk, represented in solid line, is here to be read 
from right to left and records the succession of buds and leaves clockwise around 
the tree from the root.  Starting from height $0$, it makes an ascending step for each 
encountered leaf and a descending one for each bud, hence ends at height $1$. Its depth
is here $1$, equal to the absolute value of the minimum visited height. The surrounding
walk, represented in dotted line, is to be read from left to right and records
to sequence of buds and blossom subtrees counterclockwise around the root vertex
from the root.  To check whether the tree at hand contributes to $R_n$, we let
the surrounding walk start from height $n$ and make an ascending step for
each encountered bud and a descending one for each blossom subtree. Each descending step
may be viewed as the net contribution of the subtree within the contour walk. Here each
blob corresponds to such a subtree. The depth of the contour walk will remain less
or equal to $n$ if the depth of each (vertically shifted) contour walk in a blob starting 
at height $i$ remains less or equal to $i$. }{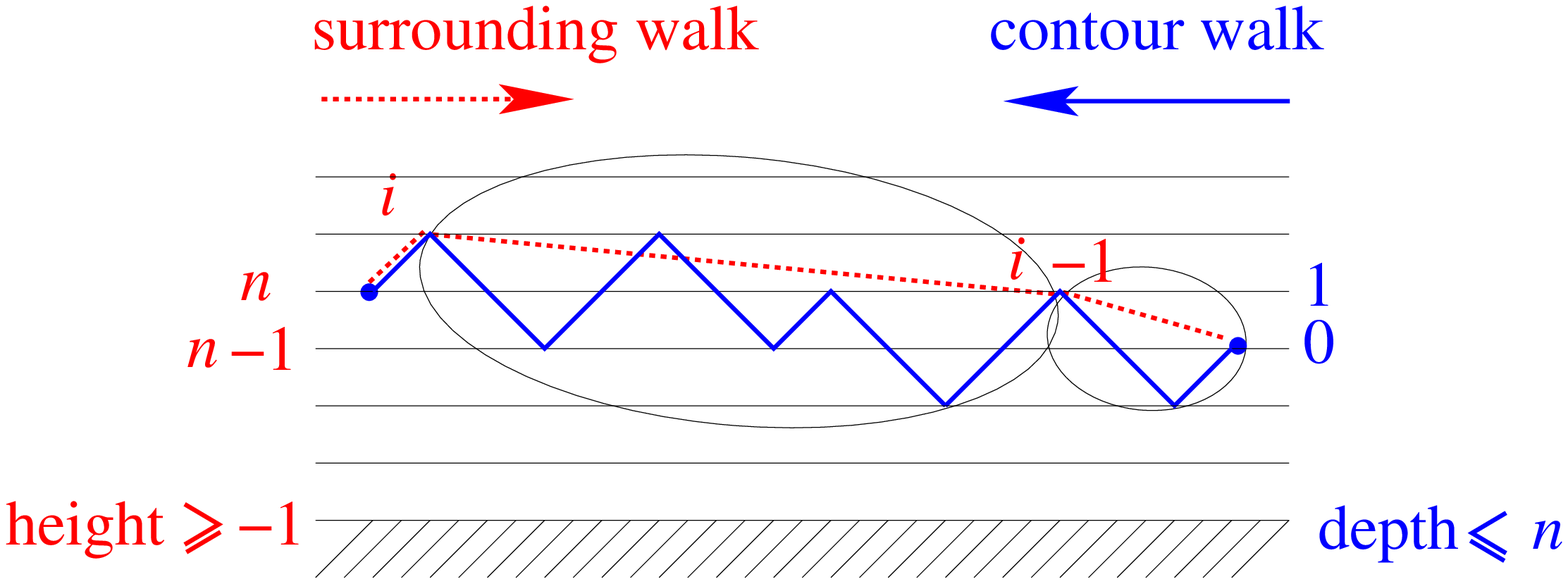}{11.cm}
\figlabel\contour
The inverse of the cutting procedure consists in reading the sequence of buds and leaves
in counterclockwise order around the tree (see Fig.\closing-(a)) and connecting each pair of a bud 
immediately followed by a leaf into an edge (see Fig.\closing-(b))  The process is repeated until only
one leaf is left (see Fig.\closing-(c)), which is chosen as the incoming leg while the root becomes
the outcoming one (see Fig.\closing-(d)). Note that this root is usually encircled by a number of 
edges separating it from the outer face. This number is nothing but the geodesic
distance between the legs. A simple way to measure it is to
forbid to encircle the root in the closing process. The geodesic distance is then given
by the number of {\it excess} buds, i.e. those buds which remain unmatched. Alternatively, 
we may define the {\it contour walk} of a blossom tree 
by reading the sequence of buds and leaves clockwise around the tree starting from the root
and performing a down (resp. up) step for each encountered bud (resp. leaf)
(see Fig.\contour).
Starting from height $0$ and making steps of $\pm 1$, the contour walk ends
at height $1$ and the geodesic distance between the legs of the corresponding
diagram is given by the {\it depth} of the walk, defined as the absolute value 
of the minimum height reached by the walk.

\subsec{Master equation}
We may now define the generating function $R_n$ for two-leg diagrams with
weight $g_k$ per $2k$-valent vertex and with legs at a geodesic distance 
{\it less or equal} to $n$. From the above bijection, $R_n$ is also the generating 
function for blossom trees with weight $g_k$ per $2k$-valent inner vertex 
and whose contour walk has depth less or equal to $n$ or equivalently with 
at most $n$ excess buds in the closing process that forbids encircling the root.

It is now easy to write down a {\it master equation} for $R_n$ by 
decomposing each blossom tree according to the environment of the inner vertex 
attached to the root (root vertex). This environment may be decomposed into a sequence of buds and 
{\it descendent subtrees}. All these subtrees are themselves blossom trees either 
made of a single leaf or satisfying the above characterizations (B1-B3). 
Note that the contour walk of the original tree is the concatenation of down
steps for the buds and of the (vertically shifted) contour walks of the blossom subtrees. Hence, the 
constraint of depth less than $n$ translates into constraints on the
depth reached by the contour walk of each blossom subtree. This turns into a closed 
relation between the $R_n$'s. 

More precisely, the environment of the root vertex is best described by yet another walk, 
the {\it surrounding walk}, now obtained by reading the sequence of buds and subtrees 
in the opposite, counterclockwise direction around the root
vertex. Starting from height $n$ at the root, we make a step $+1$ for each encountered 
bud and $-1$  for each encountered blossom subtree (see Fig.\contour). Each $-1$ step
may be seen as the net contribution of the blossom subtree at hand within the total 
(reversed) contour walk. 
To ensure that the depth of the entire contour walk remains 
less or equal to $n$, we must require that each subtree encountered at height $i$ in
the surrounding walk has a contour walk with depth less or equal to $i$.
This leads us to attach a weight $R_i$ to each $-1$ step from height $i$ to height $(i-1)$.

The description of the root vertex environment by its surrounding walk and the attached weights 
are efficiently encoded into a ``transfer matrix" $Q$ acting on a formal orthonormal basis 
$|i \rangle$ ($i\in \IZ$) as
\eqn\transfer{Q|i\rangle= |i+1\rangle + R_i\, |i-1\rangle}
The generating function $Z_{a,b}(k)$ for walks of $k$ steps going from, say height $a$ 
to height $b$, and with weight $R_i$ for each descent $i\to (i-1)$ is easily 
expressed in terms of $Q$ as
\eqn\zab{Z_{a,b}(k)\equiv \langle b | Q^k |a \rangle}
Note that $Z_{a,b}(k)$ is non zero only for $(b-a)=k\ {\rm mod}\ 2$.
{}From the above analysis, the quantity $g_{k}\, Z_{n,n-1}(2k-1)$ is the generating function 
for blossom trees with a $2k$-valent root vertex and with contour
walk of depth less or equal to $n$.
To incorporate arbitrary (even) valences, we introduce the quantities
\eqn\Vprime{V'(Q)=\sum_{k\geq 1}  g_{k} Q^{2k-1}}
and 
\eqn\Vprimeab{V'_{a,b}\equiv \langle b | V'(Q) |a \rangle= \sum_{k\geq 1} g_{k} Z_{a,b}(2k-1)}
(which is non zero only for $(b-a)$ odd). The quantity $V'_{a,b}$ may be viewed as
a ``grand-canonical" generating function for walks from $a$ to $b$ and with arbitrary
odd length.

With these notations, we may finally write the master equation
\eqn\master{\encadremath{R_n=1+V'_{n,n-1}}}
where the first contribution ($1$) corresponds to having a single
leaf connected to the root. This relation is valid for all $n\geq 0$ with
the convention that $R_{i}=0$ for all $i<0$. This equation determines
all the $R_n$'s as formal power series of the $g_k$'s upon imposing that
$R_n=1+{\cal O}(\{g_k\})$ for all $n\geq 0$.
\fig{A schematic representation of the truncated master equation (2.6). The generating
function $R_n$ is decomposed according to the environment of the root vertex,
with $0$, $1$, $2$ buds and $1$, $2$, $3$ blossom subtrees according to its valence
$2$, $4$,  $6$ respectively. To ensure that the depth of the
contour walk does not exceed $n$, we have to impose bounds on the maximal number of excess 
buds in each subtree. We have represented (here for $n=3$) these excess buds in maximal number
for each subtree. The index of $R$ for each subtree is nothing but this
maximal number.}{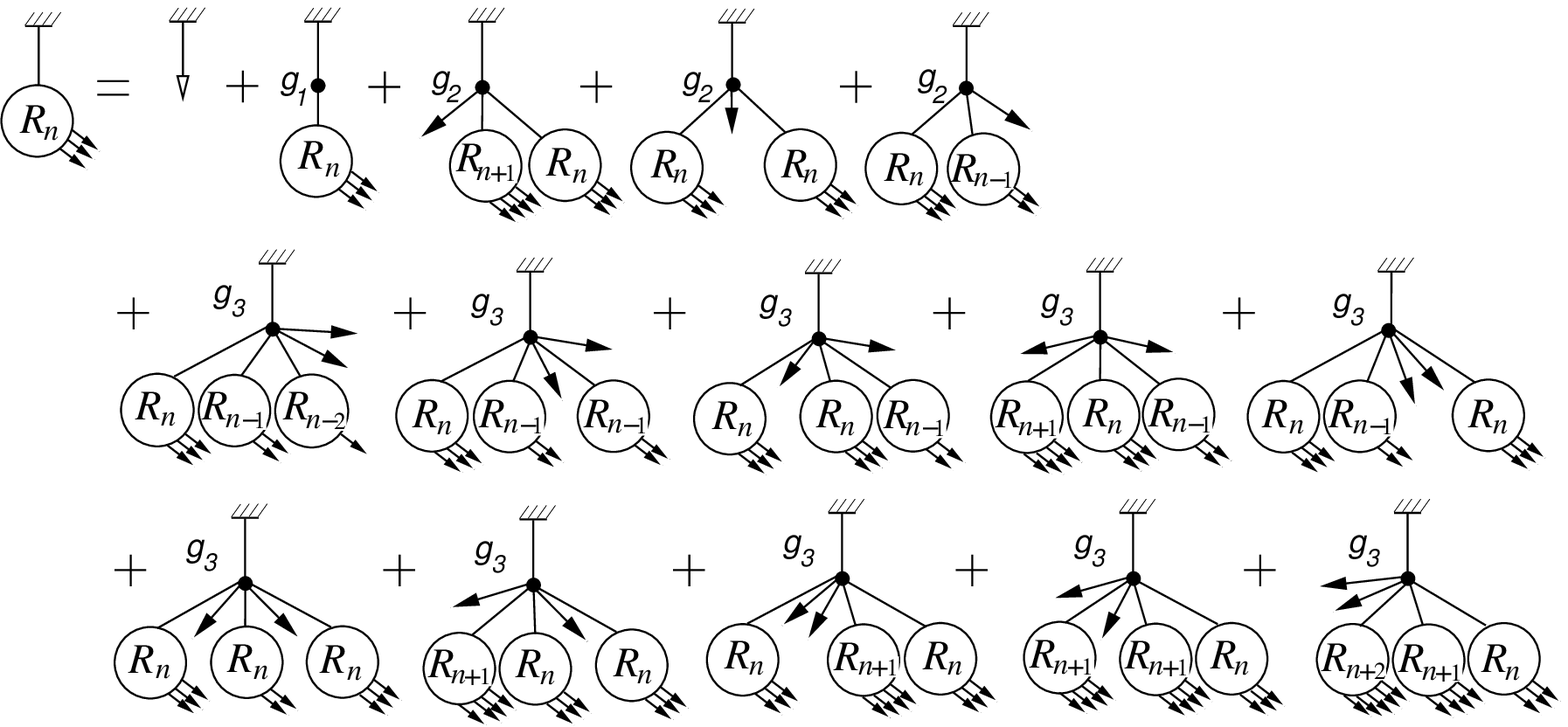}{13.5cm}
\figlabel\masterfig
\noindent For illustration, let us consider the simple case of diagram with only bi-, tetra- 
and hexavalent inner vertices, i.e. with $g_k=0$ for all $k\geq 4$.
The master equation then reads
\eqn\masterbth{\eqalign{R_n & = 1  +  g_1 R_n + g_2 R_n (R_{n-1}+R_n+R_{n+1}) \cr
& + g_3 R_n \big(R_{n-2}R_{n-1} +R_{n-1}^2 + R_{n}^2 + 2 R_n (R_{n-1}+R_{n+1})
\cr & \qquad \qquad \qquad \qquad \qquad + R_{n-1}R_{n+1} +R_{n+1}^2 + R_{n+1}R_{n+2}\big)\cr}}
This relation is illustrated in Fig.\masterfig. 
In the following, we will also make use of the generating function $R$ for two-leg
diagrams without restriction on the geodesic distance. This may be recovered
as the limit $n\to \infty$ of $R_n$, hence $R$ satisfies the relation
\eqn\eqforR{R=1+\sum_{k\geq 1} g_k {2k-1 \choose k} R^{k}}
obtained by counting walks of length $2k-1$ with global height difference equal
to $-1$, thus having $k$ descending steps each weighted by $R$.
In the truncated case above, Eq.\eqforR\ reduces to
\eqn\eqforRtrunc{R=1+g_1 R+3 g_2 R^2 +10 g_3 R^3}

\subsec{Other applications of the master equation}
The master equation may be rewritten as
\eqn\labelled{R_n={1\over 1-V'_{n-1,n}}}
by noting that $V'_{n,n-1}=R_n V'_{n-1,n}$. The latter is a direct consequence of the
general reflection relation
\eqn\reflec{Z_{a,a-1}(2k-1)=R_a Z_{a-1,a}(2k-1)}
obtained by first transferring the weights $R_i$ to the ascending steps $(i-1) \to i$
and then reflecting the walks, thus interpreting them as walks from
$(a-1)\to a$. Note that there are as many descending steps $i\to (i-1)$ as ascending steps 
$(i-1)\to i$ except for $i=a$ where there is one more descent (before reflection), hence the 
extra factor $R_a$ needed to reproduce the correct weight.

\fig{An example of labeled mobile (left) made of (inflated) polygons attached by their
vertices into a tree-like structure. The vertex labels around each polygon obey
the rule displayed on the right that they either decrease by $1$ or increase weakly
clockwise.}{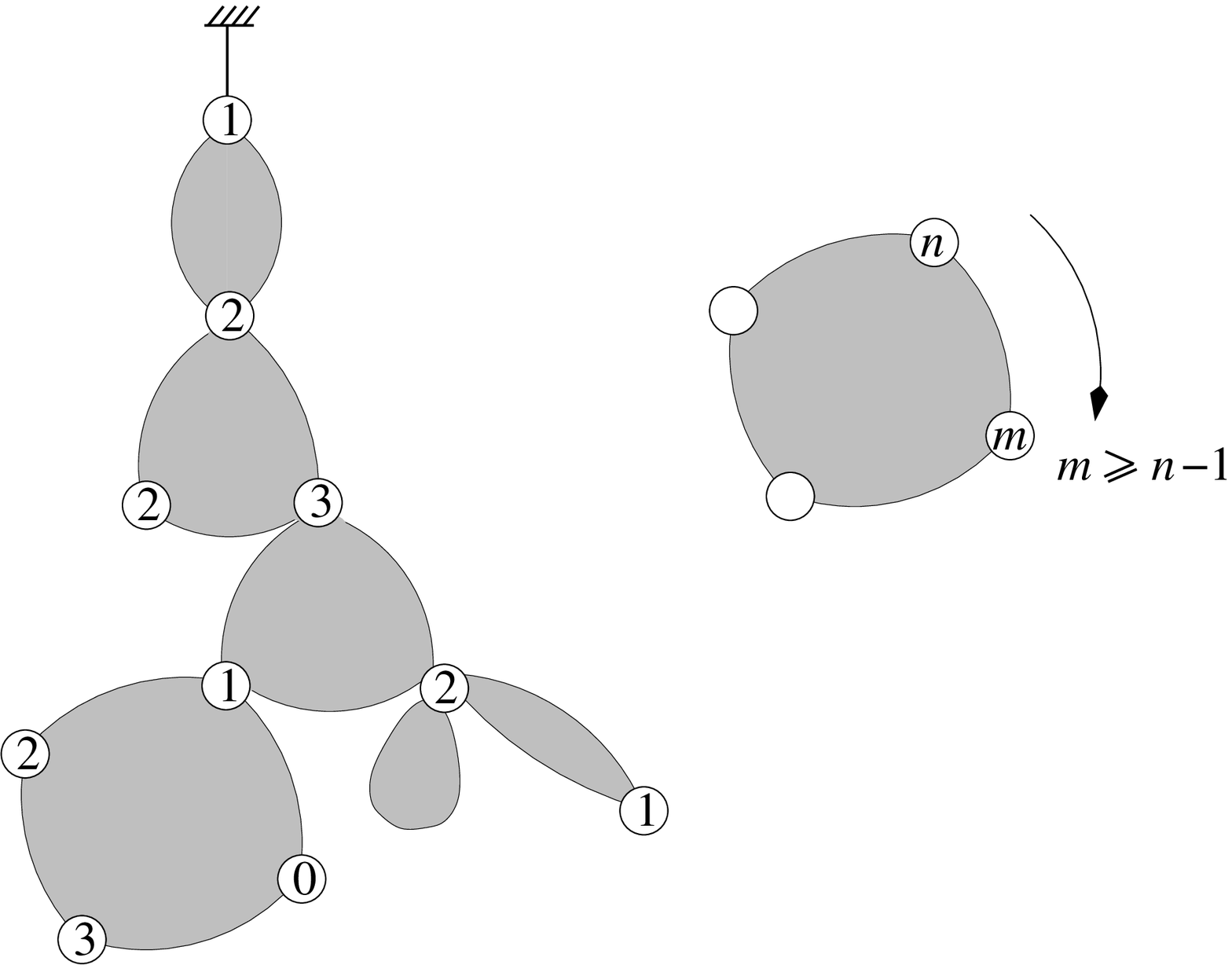}{8.cm}
\figlabel\samplemobile
In its alternative form \labelled, the master equation appeared in a slightly different 
enumeration problem, that of so-called {\it labeled mobiles} \MOB, generalizing
the so-called well labeled trees of Ref.\CORA.
These mobiles are made of rigid polygons, with weight $g_k$ per $k$-gon,
glued by their vertices into a rooted tree-like object (see Fig.\samplemobile\ for an example).
The labels are subject to local rules around each $k$-gon: going
clockwise around a $k$-gon, the labels must either decrease by $1$ or increase weakly. 
This rule again is efficiently described by the transfer matrix $Q$ of Eq.\transfer\ by reading 
the labels around $k$-gons
from the vertex at which they are suspended and translating the sequence of labels into heights 
forming a closed walk of $k$ steps around the $k$-gon. This walk is then transformed into a closed 
{\it boundary walk} of length $2k$ by transforming each non-negative step $(+p)$ into a descent followed
by $p+1$ unit ascending steps. Denoting by $R_n$ the generating function for planted mobiles
with root labeled $n$ and inspecting the possible environments of the root, we have the relation
\eqn\RtoM{R_n={1\over 1-\sum\limits_{k\geq 1} g_k M_n(k)}}
where $M_n(k)$ is the generating function for $k$-gons rooted at a vertex labelled $n$
and with $k-1$ sub-mobiles dangling at the other vertices. In $M_n(k)$, we have
to assign a weight $R_i$ to each vertex of the $k$-gon labelled $i$, except for the root
vertex. This amounts equivalently to assign a weight $R_i$ to each $-1$ step from 
height $i$ in the boundary walk, except for the $-1$ step from the root. 
This leads to
\eqn\MtoZ{M_n(k)=Z_{n-1,n}(2k-1) =\langle n|Q^{2k-1}|n-1>}
Eqs.\RtoM\ and \MtoZ\ boil down to Eq.\labelled. 
Again this equation must be supplemented with some initial conditions, 
for instance that $R_i=0$ for $i<0$ if we demand that the labels remain non-negative. 

Returning to graph enumeration problems, it was shown that mobiles with
non-negative labels are in one to one correspondence with planar graphs with {\it faces} of
even valences only and with a {\it distinguished origin vertex} [\xref\MOB,\xref\CORA]. 
The vertices of the mobile are in correspondence with vertices of the associated graph, and
the labels represent the geodesic distance (on the associated graph) of these vertices from the 
origin vertex. In this framework, $R_n$ is now interpreted as the generating function for the graphs 
with a {\it distinguished edge} at geodesic distance less than $n$ from the origin vertex.

Finally, a third domain of application of the master equation comes from the study of
spatially extended branching processes, i.e. probabilistic models for the
evolution and spreading of a population. In this language, the index $n$ stands for
the (discrete one-dimensional) position of individuals and the blossom trees or
mobiles are interpreted as the genealogical structure of families. As shown 
in Ref.\ONEWALL, the generating function $R_n$ is then related to the probability for a 
population with germ at position $n$ never to spread up to position $0$.

\subsec{Integrability of the master equation}
A remarkable feature of the master equation is that it may be solved exactly.
More precisely, we will consider truncated versions of Eq.\master\ by
considering only diagrams with valences up to, say $2m$. This amounts to
set $g_k=0$ for $k>m$ in which case $V'(Q)$ is an odd polynomial of degree $2m-1$.
As explained in Ref.\GEOD, the general solution of Eq.\master, now with 
arbitrary initial conditions but with the requirement that it converges at
$n\to \infty$ takes the surprisingly simple form
\eqn\sol{R_n=R\, {u_{n+1} u_{n+4}\over u_{n+2} u_{n+3}}}
where $u_n$ has a multi-soliton structure involving $m-1$ integration constants
$\lambda_1,\ldots,\lambda_{m-1}$.  Those constants may be fixed by the initial conditions. 
This structure is characteristic of (discrete) integrable systems and implies
the existence of conserved quantities. More precisely,  
by inverting the system $R_{i}=R_{i}(\{\lambda\})$ for any set 
$i=n,n+1,\ldots,n+m-2$, one may {\it in principle} construct $m-1$ conserved quantities for 
Eq.\master\ in the form $\lambda_j=\Lambda_j(R_n.R_{n+1},\dots,R_{n+m-2})=$ 
const.\ independently of $n$. In practice, this construction is however difficult
to implement. The remainder of the paper is devoted to the explicit construction 
of conserved quantities {\it ab initio}, i.e. without reference to the
above solution \sol. This construction will be carried out by use of simple 
combinatorial tools only.

\newsec{Random walks and heaps of dimers}
\subsec{From random walks to heaps of dimers}
In this Section, we shall investigate a number of properties of the partition
function of random walks $Z_{a,b}(k)$ above. In particular, we will 
emphasize the connection between random walks and so-called {\it heaps of dimers},
leading us eventually to an inversion relation involving partition functions 
for hard dimers on a segment. All the present relations will be instrumental for proving
the conservation statements of Section 4 below.
\fig{A schematic representation of the fundamental identity (3.1) for walks.
The walks are decomposed according to either their first step (first line)
or their last one (second line).}{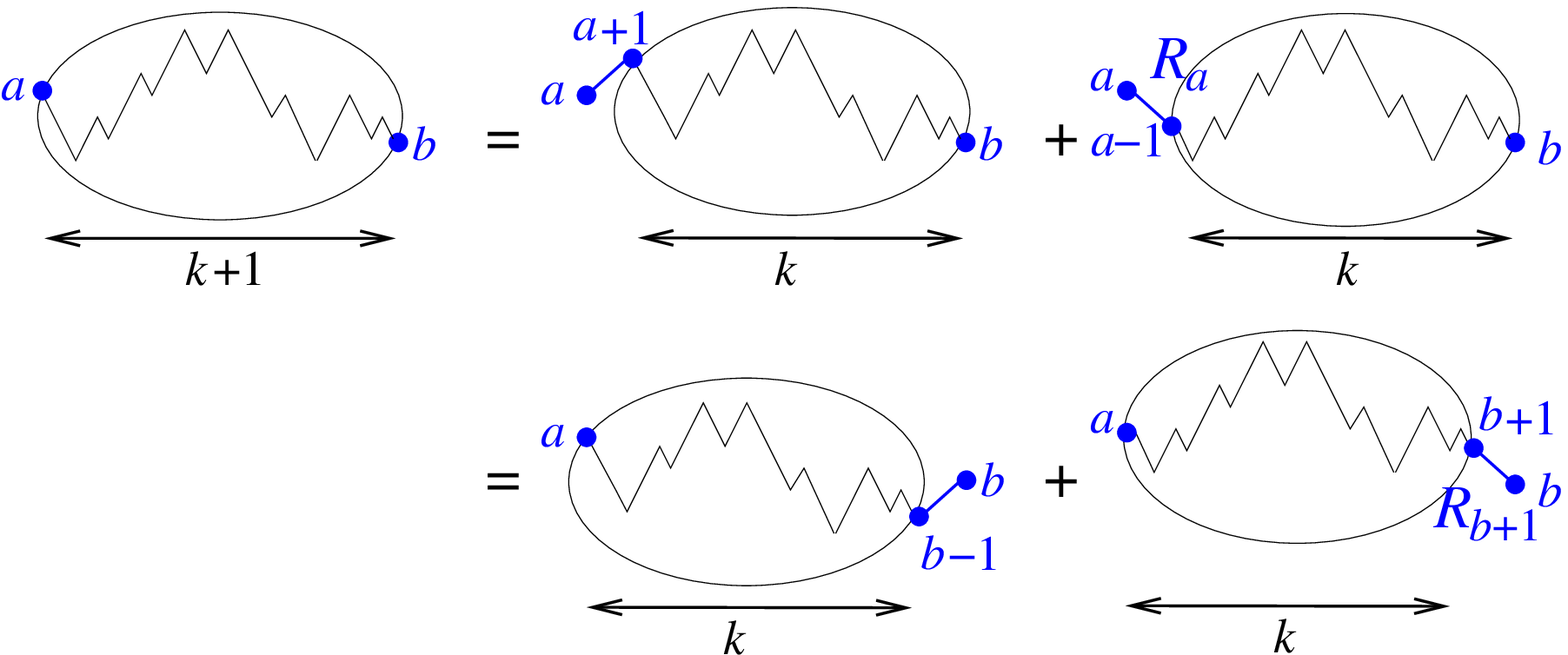}{13.cm}
\figlabel\zab
A first fundamental identity is obtained by expressing the generating function 
for paths of length $k+1$ in two different ways: (i) as the concatenation of
a path of length $1$ (first step) and a path of length $k$, or (ii) as the concatenation of
a path of length $k$ and a path of length $1$ (last step). According to whether the first
(resp. the last) step is up or down, we have respectively
\eqn\fundZ{Z_{a,b}(k+1)=Z_{a+1,b}(k)+R_a Z_{a-1,b}(k) = Z_{a,b-1}(k)+R_{b+1}Z_{a,b+1}(k)}
as illustrated in Fig.\zab.
\fig{An example of a heap of dimers (a) and its canonical representation (b).
The thickened dimers correspond to the right projection of the heap, forming
in (b) a hard dimer configuration on the vertical line $x=0$.  By rotating the picture 
(b) and extending each dimer into a $2\times 1$ rectangle, we obtain a piling up of these 
rectangles (c) whose bottom row is formed by the former right projection.}{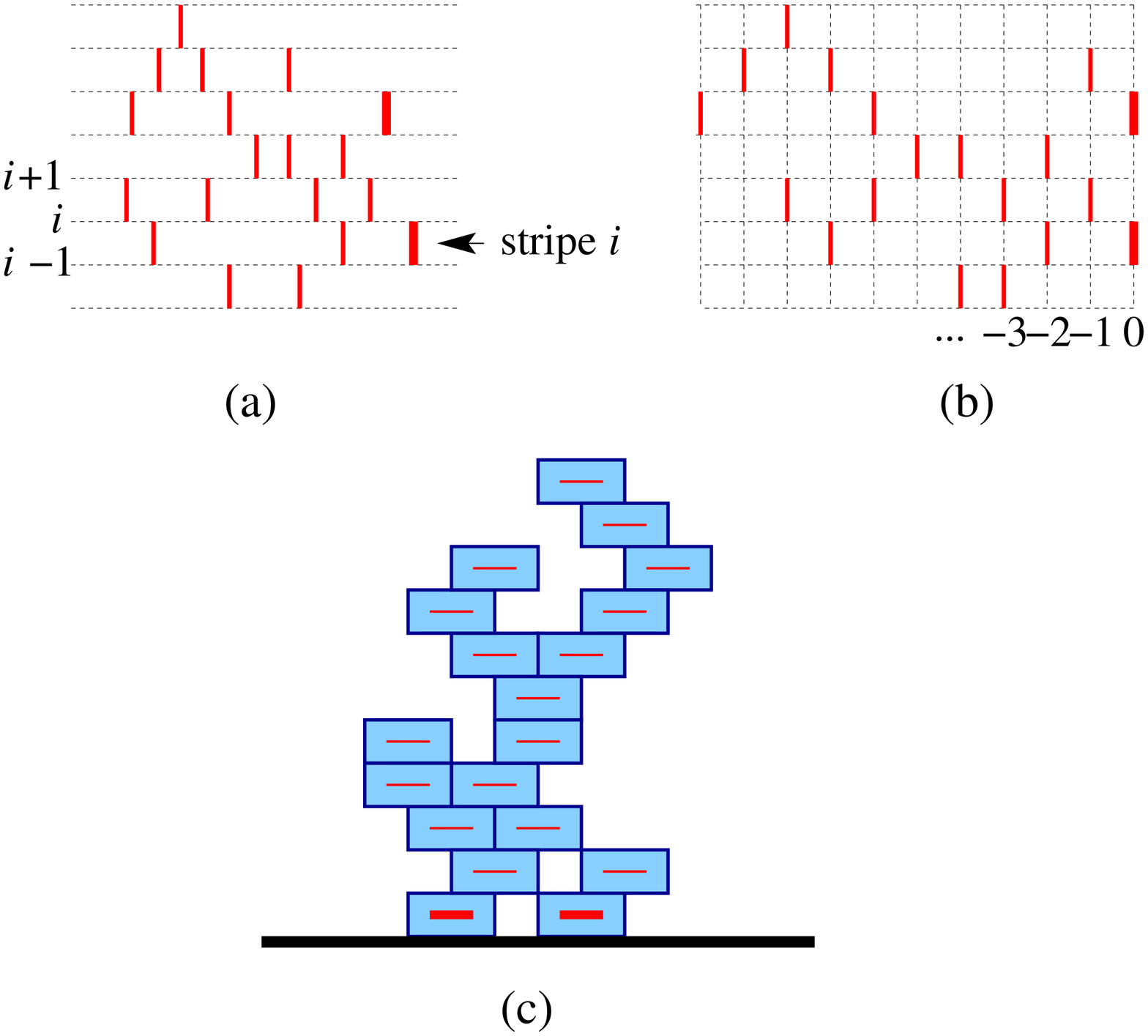}{12.cm}
\figlabel\heaps
A second important property is that we may interpret $Z_{a,b}(k)$ as a generating function 
for heaps of dimers, a particular case of the so-called heaps of pieces introduced in
Ref.\VIENNOT.  These heaps of dimers also occur in relation to lattice animals
\BMR\ and so-called Lorentzian gravity \DGL\ and may be defined as follows. Let us 
decompose the plane into parallel (say horizontal) stripes of width $1$ 
by drawing horizontal lines at integer vertical positions.  The stripes are labeled 
by the position of their top boundary. We then place within the stripes a number 
of dimers, i.e.  vertical segments of length $1$ (see Fig.\heaps-(a)).  All dimers
may freely slide within their stripe provided they do not cross a dimer within
the same stripe or within the stripe immediately above and below. A heap of dimers
is defined modulo this sliding freedom and only records the {\it relative} positions
of the dimers. In a canonical representation, we may place all the dimers 
at integer horizontal coordinates within their respective stripes with the requirement
that these coordinates be negative or zero and {\it maximal} (see Fig.\heaps-(b)). 
This is realized by pushing all dimers as much to the right as possible while staying 
at integer horizontal positions in the left half plane.
This in turn allows upon rotating the picture by $90^\circ$ clockwise to view a configuration 
as a piling up of rectangles of size $2\times 1$ (see Fig.\heaps-(c)) deposited on the
$x=0$ line.
For each heap of dimers, we may define its {\it right projection} as the configuration
made of those dimers which may slide freely all the way to infinity to the right. 
In the canonical representation, these are those dimers with horizontal position $0$.
Clearly, the right projection of any heap of dimers is a configuration of {\it hard dimers}
on a (vertical) line, i.e. a set of dimers occupying (vertical) segments $[i-1,i]$
with the restriction that no two dimers may come in contact (the occupied segments must be 
disjoint). 
\fig{An example of a pyramid (a) and a half-pyramid (b), both with base $[a,b]$.}{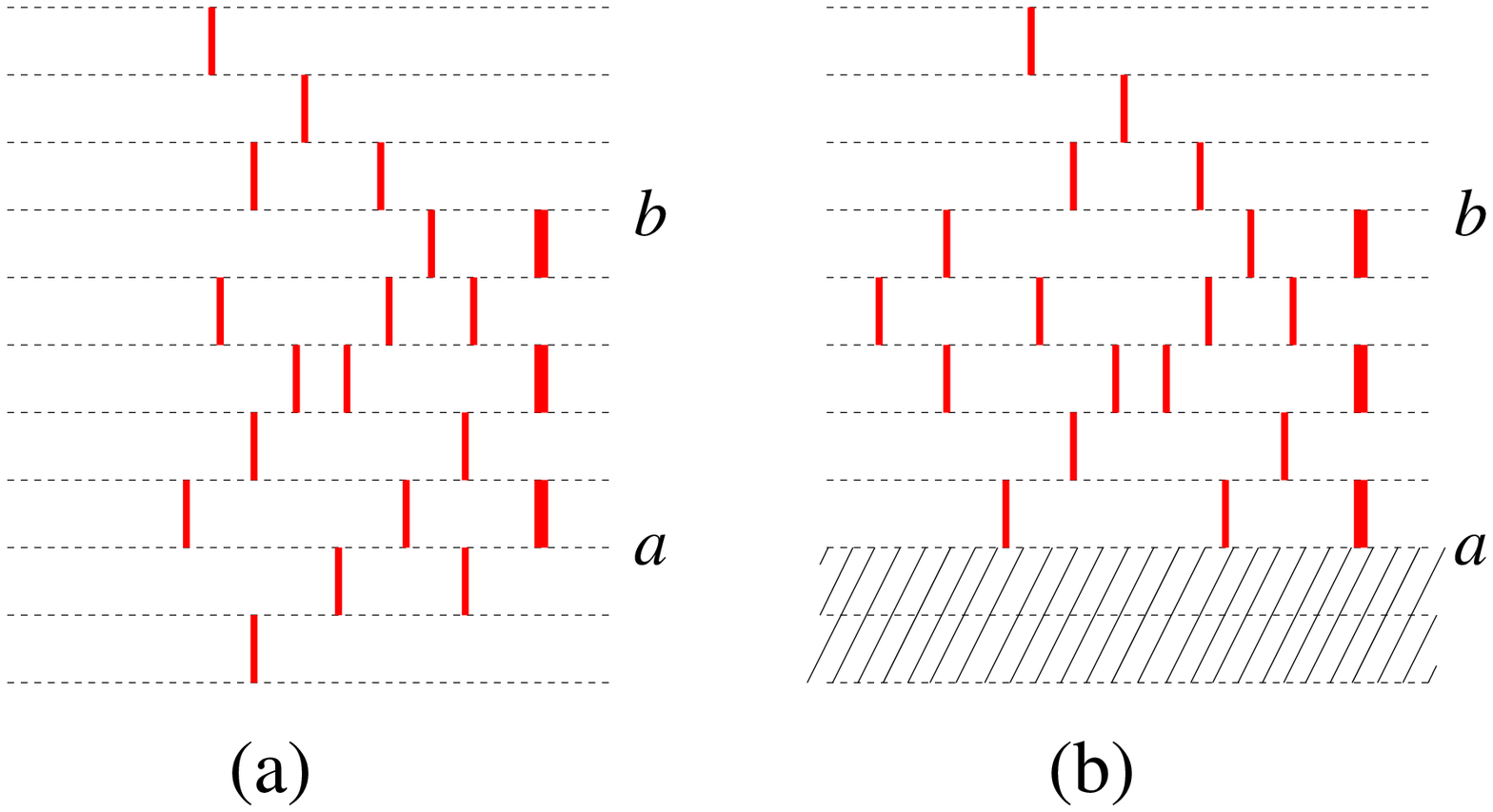}{10.cm}
\figlabel\pyramid
For $(b-a)$ odd and positive, we may consider the maximally
occupied hard dimer configuration of the segment $[a,b]$. It is made of $(b-a+1)/2$ 
dimers occupying the elementary segments $[a+2i,a+2i+1]$, $i=0,\ldots,(b-a-1)/2$. 
Any heap of dimer having this maximally occupied segment as right projection will be
called a {\it pyramid} of base $[a,b]$ (see Fig.\pyramid-(a) for an example). We then define the 
generating function $H_{a,b}(k)$ for pyramids
of base $[a,b]$ made of a {\it total} number $k$ of dimers, with weight $R_i$ per dimer
in the stripe labeled $i$, {\it except} for the dimers of the right projection. Note that
$H_{a,b}(k)$ is non zero only if $k\geq (b-a+1)/2$ and that $H_{a,b}((b-a+1)/2)=1$.
\fig{A schematic representation (first line) of the equality (3.2) between generating functions
$Z_{a,b}(2k-1)$ for walks from $a$ to $b$ with $2k-1$ steps and $H_{a,b}(k)$ for
pyramids of base $[a,b]$ with a total of $k$ dimers. The similar
representation (second line) for the relation (3.3) between the generating function 
$Z^{+}_{a,b}(2k-1)$ of positive walks and that $H^{+}_{a,b}((k)$ of 
half-pyramids.}{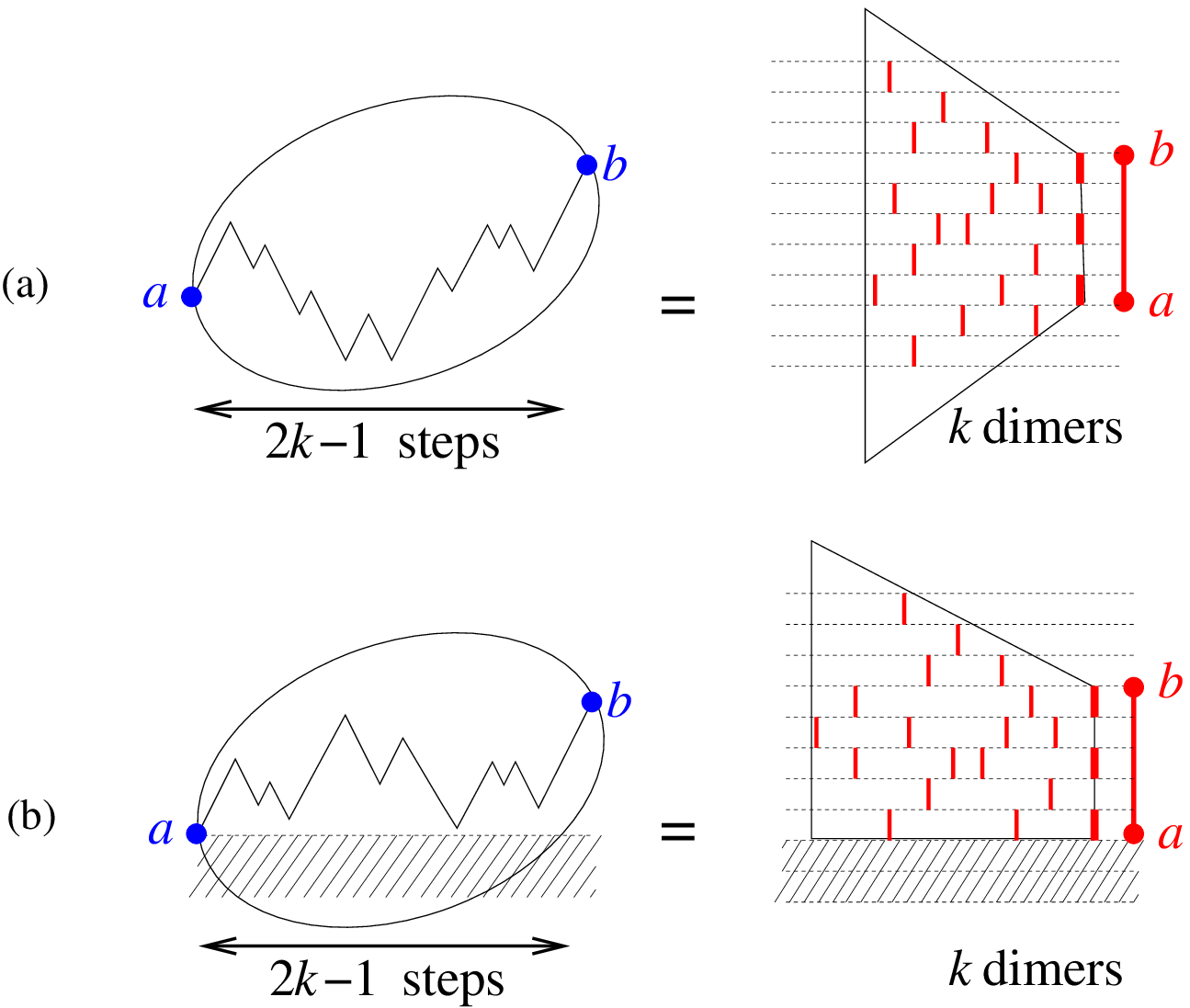}{12.cm}
\figlabel\rwtopy
We these definitions, we have the remarkable relation
\eqn\ZtoH{Z_{a,b}(2k-1)=H_{a,b}(k)}
valid for $(b-a)$ odd and positive, and $k\geq 1$. This relation is illustrated in
Fig.\rwtopy-(a). It follows from a general bijection between walks and heaps of pieces
given in Ref.\VIENNOT. This bijection is recalled in detail in Appendix A below 
for the particular case at hand.

In the next section, we shall also make use of the generating function $Z^{+}_{a,b}(k)$ for walks of 
$k$ steps going from height $a$ to height $b$, with weight $R_i$ per descent $i\to (i-1)$, and 
{\it which stay at or above height $a$}. In the following, we shall refer to these walks as
{\it positive walks}. The generating function $Z^{+}_{a,b}(k)$ may be obtained from
$Z_{a,b}(k)$ by simply taking $R_i\to 0$ for all labels $i\leq a$. In particular, we deduce that,
for $(b-a)$ odd and positive and $k\geq 1$
\eqn\ZtoHplus{Z^{+}_{a,b}(2k-1)=H^{+}_{a,b}(k)}
where $H^{+}_{a,b}(k)$ denotes the generating function for {\it half-pyramids} of
base $[a,b]$ with a total of $k$ dimers, i.e. pyramids with dimers only in stripes
with labels $i>a$ (see Fig.\pyramid-(b) for an example). This relation \ZtoHplus\ is illustrated
in Fig.\rwtopy-(b).

\subsec{Inversion relation}

Viewing hard dimers on a line as mutually excluding and fermionic objects (with at most
one dimer in a unit segment), heaps of dimers appear as their interacting bosonic counterpart 
in which arbitrarily many dimers may be piled up in each stripe. It is therefore natural 
to expect some boson/fermion inversion relations to hold between their respective generating
functions. An example of such inversion relation in a ``grand-canonical" ensemble,
i.e. for arbitrarily many dimers, may be found in Refs.[\xref\VIENNOT,\xref\DGL].
In this Section, we shall derive another inversion formula, now in the canonical ensemble,
i.e. with fixed numbers of dimers. 

Introducing the generating function  $\Pi_{a,b}(k)$ for {\it configurations of $k$ hard dimers} 
in the segment $[a,b]$ and with weight $R_i$ per dimer in the segment $[i-1,i]$, 
we have the remarkable inversion relation 
\eqn\geninversion{\sum_{\ell=0}^{k}(-1)^{k-\ell} \Pi_{a,a+2 k-1}(k-\ell)Z^{+}_{a,a+2j}(2\ell) =\delta_{k,j}} 
for $k,j\geq 0$, with the convention that $\Pi_{a,a-1}(0)=1$.  
This relation may be written in a more compact matrix form by defining a lower-triangular, 
semi-infinite matrix ${\bf Z}(a)$ with entries
\eqn\matZ{{\bf Z}(a)_{i,j}\equiv Z^{+}_{a,a+2j}(2i)} 
for $0\leq j\leq i$, and a lower-triangular, semi-infinite matrix ${\bf D}(a)$ with entries
\eqn\matD{{\bf D}(a)_{k,i}=(-1)^{k-i}\Pi_{a,a+2k-1}(k-i)}
for $0\leq i\leq k$ while all other entries vanish. The relation \geninversion\ now reads \eqn\matinv{{\bf D}(a){\bf Z}(a)={\bf I}}
with ${\bf I}$ the (semi-infinite) identity matrix. To prove the relation \geninversion, 
let us consider {\it pairs} ${\cal P}$ made of
\item{(i)} a hard dimer configuration in the segment $[a,a+2k-1]$
\item{(ii)} a half-pyramid with projection $[a,a+2j-1]$
\par
\noindent with a total number of $k$ dimers. In the half-pyramid, each dimer
in the stripe $i$ (including the dimers of the right projection) receives a weight $R_i$.
On the contrary, in the hard dimer configuration, each dimer in the 
segment $[i-1,i]$ receives a weight $-R_i$. Denoting by $\ell$ the number of dimers in the half-pyramid,
the generating function $P_{a}(j,k)$ for the above pairs reads
\fig{An example (left) of pair ${\cal P}$ of a hard-dimer configuration
in the segment $[a,a+2k-1]$ and a half-pyramid of base $[a,a+2j-1]$ (with $j\leq k$),
with a total number $k$ of dimers (here $k=7$). These are concatenated (right) 
into a larger heap ${\cal H}$. The equivalence class ${\cal C}({\cal H})$ of those
pairs leading to the same heap ${\cal H}$ is constructed by considering in ${\cal H}$ the dimers
which belong to the left projection of the heap but do not belong to the right
projection. In the present example, there is exactly one such (encircled) dimer. The
equivalence class ${\cal C}({\cal H})$ is generated by distributing in all possible ways these dimers
either in the hard-dimer configuration or in the half-pyramid. If we assign opposite weights to dimers
in the hard-dimer configuration and in the half-pyramid, this results
in a vanishing net contribution of the class ${\cal C}({\cal H})$ unless the 
left and right projections of ${\cal H}$ are identical. This happens only when ${\cal H}$ is 
the maximally occupied hard-dimer configuration made of $k$ dimers in the segment 
$[a,a+2k-1]$.}{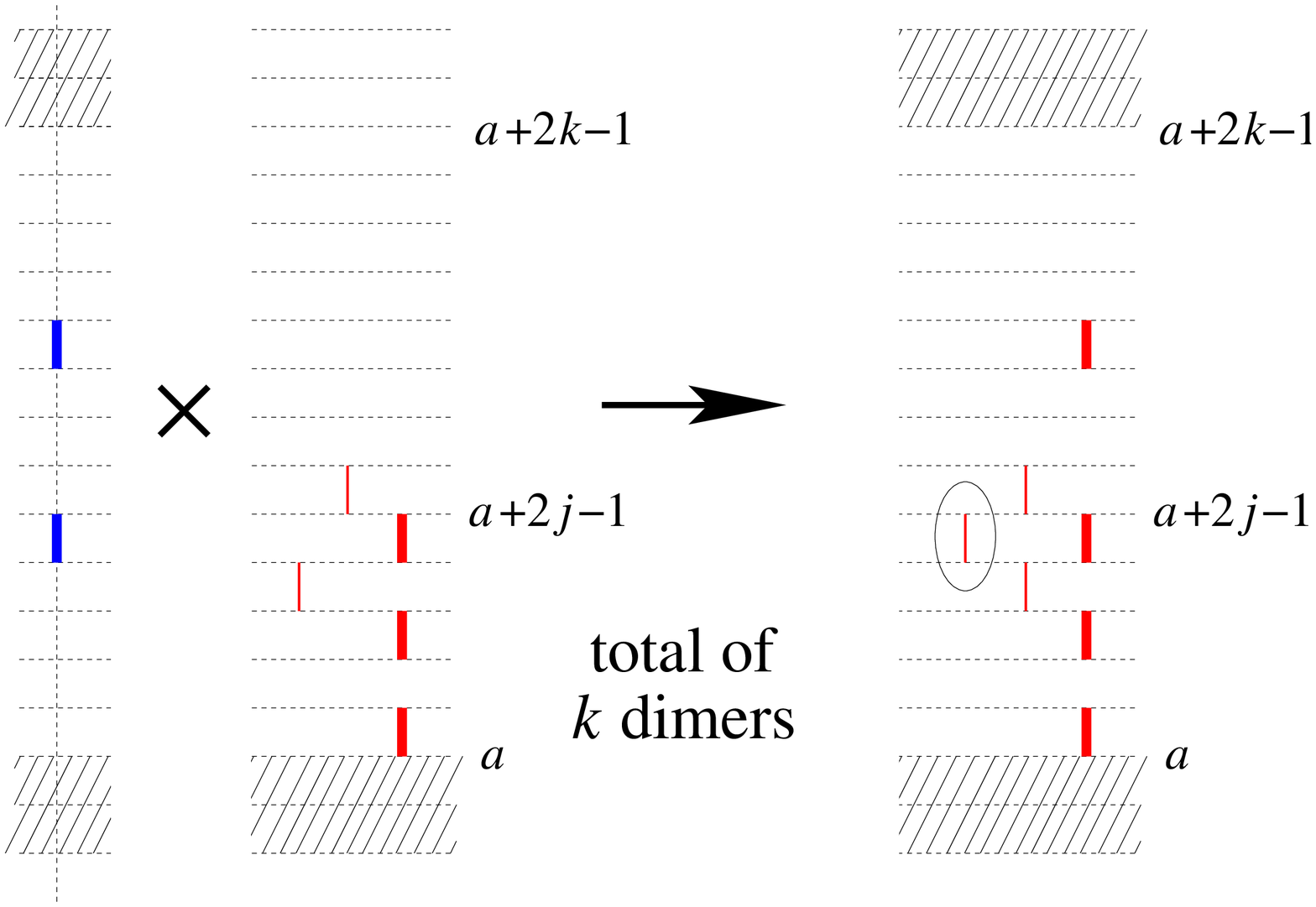}{8.cm}
\figlabel\pairs
\eqn\generpair{P_{a}(j,k)=\Pi_{a,a+2j-1}(j) \sum_{\ell=0}^{k}(-1)^{k-\ell} \Pi_{a,a+2 k-1}(k-\ell)
H^{+}_{a,a+2j-1}(\ell) }
Note that $P_{a}(j,k)$ is non zero only if $k\geq j$.
For each such pair ${\cal P}$, we may absorb all dimers of the hard dimer configuration into a larger
heap ${\cal H}({\cal P})$ of $k$ dimers by simply adding them to the left of the original 
half-pyramid (see Fig.\pairs). 
We can then regroup all pairs leading to the same larger heap ${\cal H}$ into an equivalence class
${\cal C}({\cal H})$.  Conversely, given ${\cal H}={\cal H}({\cal P})$, we may reconstruct all pairs 
in the class ${\cal C}({\cal H})$ by considering the set of dimers {\it which belong to the left projection} 
of ${\cal H}$ 
(those dimers which may be pushed all the way to infinity to the left) but {\it which do not belong to 
the right projection}. If this set contains at least one dimer, this dimer may be incorporated either 
in the hard dimer configuration of the pair or in its half-pyramid, contributing with
opposite weights. This causes the contribution to $P_{a}(j,k)$ of the class 
${\cal C}({\cal H})$ to vanish. The only case where no such vanishing
takes place is when the left projection of ${\cal H}$ is identical to its
right projection, in which case ${\cal H}$ itself is a hard dimer configuration
with $k$ dimers on the segment $[a,a+2k-1]$, therefore the unique maximally occupied 
hard dimer configuration, with contribution $(-1)^{k-j} \Pi_{a,a+2k-1}(k)$ 
as the $j$ lower dimers with label less than $a+2j-1$
belong to the half-pyramid and the $(k-j)$ upper ones belong to the hard dimer configuration.
We immediately deduce that 
\eqn\Pinv{P_{a}(j,k)=(-1)^{k-j} \Pi_{a,a+2k-1}(k)\ \delta_{k\geq j}}
with $ \delta_{k\geq j}=1$ for $k\geq j$ and zero otherwise.
Using $\delta_{k,j}=\delta_{k\geq j}-\delta_{k\geq j+1}$, we deduce
\eqn\generimpair{\eqalign{\delta_{k,j} \ \Pi_{a,a+2k-1}(k) & = P_{a}(j,k)+P_{a}(j+1,k) \cr
&=\sum_{\ell=0}^{k}(-1)^{k-\ell} \Pi_{a,a+2k-1}(k-\ell) \cr
& \ \ \times\left(H^{+}_{a,a+2j-1}(\ell) \Pi_{a,a+2j-1}(j) 
+ H^{+}_{a,a+2j+1}(\ell) \Pi_{a,a+2j+1}(j+1)\right)\cr}}
Using $H^{+}_{a,a+2j-1}(\ell)=Z^{+}_{a,a+2j-1}(2\ell -1)$, $H^{+}_{a,a+2j+1}(\ell)=Z^{+}_{a,a+2j+1}(2\ell -1)$
from Eq.\ZtoHplus\ and the relation $\Pi_{a,a+2j+1}(j+1)=R_{a+2j+1}\Pi_{a,a+2j-1}(j)$ for maximally
occupied segments, we obtain
\eqn\invfinal{\sum_{\ell=0}^{k}(-1)^{k-\ell} \Pi_{a,a+2k-1}(k-\ell) \left( Z^{+}_{a,a+2j-1}(2\ell -1) 
+ R_{a+2j+1} Z^{+}_{a,a+2j+1}(2\ell -1)\right) = \delta_{k,j}}
which, upon using   
\eqn\fundZpart{Z^{+}_{a,a+2j-1}(2\ell -1) + R_{a+2j+1} Z^{+}_{a,a+2j+1}(2\ell -1)= Z^{+}_{a,a+2j}(2\ell)}
from Eq.\fundZ, reduces to the desired inversion relation \geninversion.

Note that, as ${\bf D}(a)$ and ${\bf Z}(a)$ are lower-triangular, we may decide to 
truncate them to $m\times m$ matrices ${\bf D}_m(a)$ and ${\bf Z}_m(a)$ with indices 
$i$ strictly less than $m$. These matrices clearly satisfy the inversion relation
${\bf D}_m(a){\bf Z}_m(a)=I_m$, with $I_m$ the $m\times m$ identity matrix. 
For illustration, taking $m=3$, we have 
\eqn\sampleinv{\eqalign{{\bf D}_3(a)& =\pmatrix{1&0&0\cr -R_{a+1}&1&0\cr R_{a+1}R_{a+3}&
-(R_{a+1}+R_{a+2}+R_{a+3})&1\cr} \cr
{\bf Z}_3(a)& =\pmatrix{1&0&0\cr R_{a+1}&1&0\cr R_{a+1}(R_{a+1}+R_{a+2})&
(R_{a+1}+R_{a+2}+R_{a+3})&1\cr}\cr}} 
which are inverse of one-another. 

\newsec{Conserved quantities}
\subsec{Definition of the conserved quantities}
We now present compact expressions for a particular ``basis" of conserved quantities of 
the master equation \master. More precisely, we shall define below a set of
quantities $\{\Gamma_{2i}(n)\}$ for $i\geq 1$ and $n\geq 0$ satisfying 
$\Gamma_{2i}(n)$= const. independently of $n$. Of course, any functions of the
$\Gamma_{2i}$'s  are conserved as well and the precise choice of basis below
has been done in regard of its combinatorial nature. Other, alternative choices will be
discussed in Sects.6.1 and 6.2 below. 

Using the generating function $Z^{+}_{a,b}(k)$ of previous section for positive walks
and the grand-canonical generating function $V'_{a,b}$ for walks of odd length, 
we define the quantity $\Gamma_{2i}(n)$ by
\eqn\defgamma{\encadremath{\Gamma_{2i}(n)=Z^{+}_{n-1,n-1}(2i)\Gamma_0(n)
-\sum_{j=1}^{i} Z^{+}_{n-1,n-1+2j}(2i) V'_{n+2j-1,n-2}} }
for $i\geq 1$ and $n\geq 0$, with $\Gamma_0(n)$ given by
\eqn\defgammazero{\Gamma_0(n)= R_{n-1}-V'_{n-1,n-2}+\delta_{n,0}}
for $n\geq 0$. The quantity $\Gamma_0(n)$ itself is not {\it stricto sensu} a 
conserved quantity. However, the master equation precisely ensures that
$\Gamma_0(n)=1$ for $n\geq 1$ while $\Gamma_0(0)=1$ by definition 
(since $R_{-1}=V'_{-1,-2}=0$), hence the conservation of $\Gamma_0$ is a tautology.
This in turn allows to simplify the expression \defgamma\ for $\Gamma_{2i}(n)$
by substituting $\Gamma_0(n)=1$. The proof of the conservation of $\Gamma_{2i}(n)$
however is made much simpler by keeping as such the slightly more involved definition
\defgamma. 
Introducing the vectors ${\vec \Gamma}(n)$ and ${\vec V'}(n)$ with components
\eqn\vectors{{\vec \Gamma}(n)_i=\Gamma_{2i}(n)\quad {\rm and}\quad 
{\vec V'}(n)_j= (R_{n-1}+\delta_{n,0})\delta_{j,0}-V'_{n+2j-1,n-2}}
for $i,j \geq 0$, the above relations \defgamma\ and \defgammazero\ read simply
\eqn\matrixform{{\vec \Gamma}(n)={\bf Z}(n-1)\ {\vec V'}(n)}
with ${\bf Z}(a)$ defined as in \matZ. 
It will prove useful to invert this relation. This is readily performed by use of the
inversion relation \geninversion\ or \matinv, with the result
\eqn\VofD{{\vec V'}(n)= {\bf D}(n-1) {\vec \Gamma}(n)}
with ${\bf D}(a)$ defined in \matD, namely
\eqn\Vofgamma{V'_{n+2j-1,n-2}=\sum_{i=0}^{j} (-1)^{i-1} \Pi_{n-1,n+2j-2}(i) \Gamma_{2j-2i}(n) 
+\delta_{j,0} (R_{n-1}+\delta_{n,0})}
with again the convention that $\Pi_{n-1,n-2}(0)=1$. 

\subsec{Proof of the conservation}
We are now ready to prove that the quantities $\Gamma_{2i}(n)$ of Eq.\defgamma\ are 
conserved quantities of the master equation \master.
To this end, we shall now use the fundamental equation for $V'$
\eqn\fundV{V'_{a+1,b}+R_a V'_{a-1,b}=V'_{a,b-1}+R_{b+1}V'_{a,b+1}}
directly inherited from the fundamental equation \fundZ\ for $Z_{a,b}(k)$.
Taking $a=n+2j-1$ and $b=n-1$, we get the identity
\eqn\Vrule{\eqalign{V'_{(n+1)+2j-1,(n+1)-2}+R_{n+2j-1} & V'_{(n+1)+2(j-1)-1,(n+1)-2}
\cr &=V'_{n+2j-1,n-2}+R_{n} V'_{(n+2)+2(j-1)-1,(n+2)-2}\cr}}
For $j\geq 1$ and $n\geq 0$, we may substitute Eq.\Vofgamma\ into Eq.\Vrule, leading to
\eqn\Gammarule{\eqalign{\big(-\Gamma_{2j}(n+1)&+\sum_{i=1}^{j} (-1)^{i-1} \Pi_{n,(n+1)+2j-2}(i) 
\Gamma_{2j-2i}(n+1)\big)\cr
+R_{n+2j-1} & \big(\sum_{i=0}^{j-1} (-1)^{i-1} \Pi_{n,(n+1)+2(j-1)-2}(i) \Gamma_{2(j-1)-2i}(n+1) 
+\delta_{j-1,0} R_{n}\big)\cr
= \big(-\Gamma_{2j}(n) & +\sum_{i=1}^{j} (-1)^{i-1} \Pi_{n-1,n+2j-2}(i) \Gamma_{2j-2i}(n) \big)\cr
+R_{n} & \big(\sum_{i=0}^{j-1} (-1)^{i-1} \Pi_{n+1,(n+2)+2(j-1)-2}(i) \Gamma_{2(j-1)-2i}(n+2) 
+\delta_{j-1,0} R_{n+1}\big)\cr}}
where we have used $\Pi_{n-1,n+2j-2}(0)=\Pi_{n,n+2j-1}(0)=1$ for all $j\geq 0$ 
(unique empty dimer configuration). Noting that the boundary terms on both sides cancel, we
now change $i\to (i-1)$ in the last sum of each side and rearrange the factors of $\Gamma$, 
leading us to
\eqn\Gamrulbis{\eqalign{ 
\Gamma_{2j}(n+1) & -\Gamma_{2j}(n) \cr 
& = \sum_{i=1}^{j} (-1)^{i-1} \big\{ \Pi_{n,n+2j-1}(i) -R_{n+2j-1} \Pi_{n,n+2j-3}(i-1) \big\} 
\Gamma_{2j-2i}(n+1)\cr
& \ \ \ -\sum_{i=1}^{j} (-1)^{i-1} \big\{ \Pi_{n-1,n+2j-2}(i) - R_{n} \Pi_{n+1,n+2j-2}(i-1)\big\} 
\Gamma_{2j-2i}(n) \cr
& \ \ \ -\sum_{i=1}^{j} (-1)^{i} R_{n} \Pi_{n+1,n+2j-2}(i-1) \big(\Gamma_{2j-2i}(n+2)-\Gamma_{2j-2i}(n) \big)\cr
&  = \sum_{i=1}^{j} (-1)^{i-1} \Pi_{n,n+2j-2}(i) \big(\Gamma_{2j-2i}(n+1)-\Gamma_{2j-2i}(n) \big)\cr
& \ \ \ -\sum_{i=1}^{j} (-1)^{i} R_{n} \Pi_{n+1,n+2j-2}(i-1) \big(\Gamma_{2j-2i}(n+2)-\Gamma_{2j-2i}(n) \big)\cr
} }
\fig{A schematic representation of the fundamental identity (4.12). A configuration of hard
dimers on the segment $[a,b]$ may be viewed either as a configuration on
the segment $[a-1,b]$ where the unit segment $[a-1,a]$ is empty (hence the
subtraction) or as a configuration on the segment $[a,b+1]$ where the unit segment $[b,b+1]$ 
is empty.}{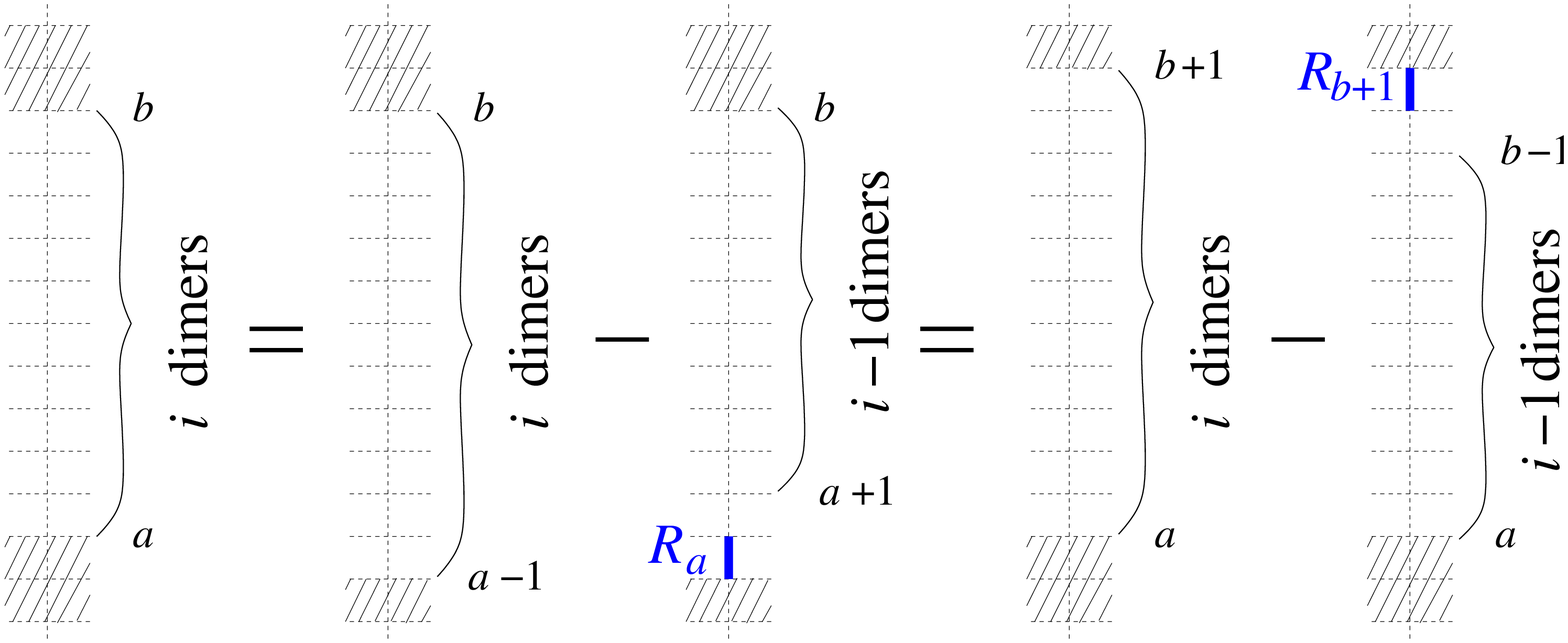}{12.cm}
\figlabel\hab
\noindent In the last equality, we have used the property
\eqn\equalhard{\eqalign{\Pi_{n,n+2j-1}(i) -& R_{n+2j-1} \Pi_{n,n+2j-3}(i-1)\cr 
& = \Pi_{n-1,n+2j-2}(i) - R_{n} \Pi_{n+1,n+2j-2}(i-1) = \Pi_{n,n+2j-2}(i)\cr} }
which is a particular instance ($a=n$, $b=n+2j-2$) of the fundamental identity
\eqn\eqhardgen{\Pi_{a,b}(i)=\Pi_{a-1,b}(i)-R_a \Pi_{a+1,b}(i-1)=\Pi_{a,b+1}(i)-R_{b+1}\Pi_{a,b-1}(i-1)}
satisfied by generating functions for hard dimers. This identity is illustrated in Fig.\hab\ and
may be viewed as the hard dimer counterpart of Eq.\fundZ\ for random walks.

{}From Eq.\Gamrulbis, we immediately deduce by induction on $j$ that
\eqn\conserv{\Gamma_{2j}(n+1)=\Gamma_{2j}(n)}
for all $j\geq 1$ and $n\geq 0$ provided that $\Gamma_0(n+1)=\Gamma_0(n)$ for all
$n\geq0$. As already mentioned, this last requirement is precisely 
guaranteed by the master equation \master. The $\Gamma_{2i}$'s therefore form
a set of conserved quantities for this equation.

For illustration, we list below the first two conserved quantities for
the truncated case of graphs with valences up to $2m=6$, corresponding to the 
truncated master equation \masterbth:
\eqn\firstcons{\eqalign{\Gamma_2(n)&= R_n -V'_{n+1,n-2}\cr
\Gamma_4(n)& = R_n(R_n+R_{n+1})-(R_n+R_{n+1}+R_{n+2})V'_{n+1,n-2}-V'_{n+3,n-2}\cr}}
with
\eqn\firstvprime{
\eqalign{ V'_{n+1,n-2}& = R_{n+1}R_n R_{n-1}\big(g_2 +g_3 (R_{n+2}+R_{n+1}+R_n+R_{n-1}+R_{n-2})\big)\cr
V'_{n+3,n-2}& = g_3 R_{n+3} R_{n+2} R_{n+1}R_n R_{n-1}\cr}}
where we explicitly substituted $\Gamma_0(n)=1$ into Eq.\defgamma. 

\newsec{Graph interpretation}
\subsec{Conserved quantities as multi-point correlation functions}
The constant value of the conserved quantities has a nice combinatorial interpretation
as a multi-point correlation function. More precisely, let us define by 
$G_{2i}(\{g_k\})$ the so-called (disconnected) $2i$-point function, i.e.
the generating function for possibly disconnected $2i$-leg diagrams, namely graphs
with inner vertices of even valences (weighted $g_k$ per $2k$-valent vertex) and 
with $2i$ legs, i.e. univalent vertices, adjacent to the same (external) face. These legs
are distinguished and labeled $1,2,\ldots,2i$ counterclockwise.
We have the identification
\eqn\valinteg{\Gamma_{2i}(n)=G_{2i}}
for all $n\geq 0$ and all $i\geq 0$ with the convention that $G_0=1$.
To prove this, we simply evaluate $\Gamma_{2i}(n)$ at $n=0$.
Noting that $V'_{2j-1,-2}=0$ as it is proportional to $R_{-1}=0$, we are left
with 
\eqn\gamzer{\Gamma_{2i}(0)=Z^{+}_{-1,-1}(2i)=Z_{-1,-1}(2i)=Z_{0,-1}(2i-1)}
where we first note that $R_{-1}=0$ automatically selects positive walks and we then remove 
the first (ascending) step with weight $1$. 
\fig{Bijection between (a) two-leg diagrams with both legs in the same face and whose outcoming 
leg is adjacent to a $2i$-valent vertex and (c) $2i$-leg diagrams. In the intermediate step
(b), we have glued the two legs into a rooted edge and cut out the $2i$-valent vertex.
The legs are labeled counterclockwise by $1,2,\ldots 2i$ with the former incoming leg 
labeled $1$. Note that the $2i$-leg diagrams need not be connected in general.}{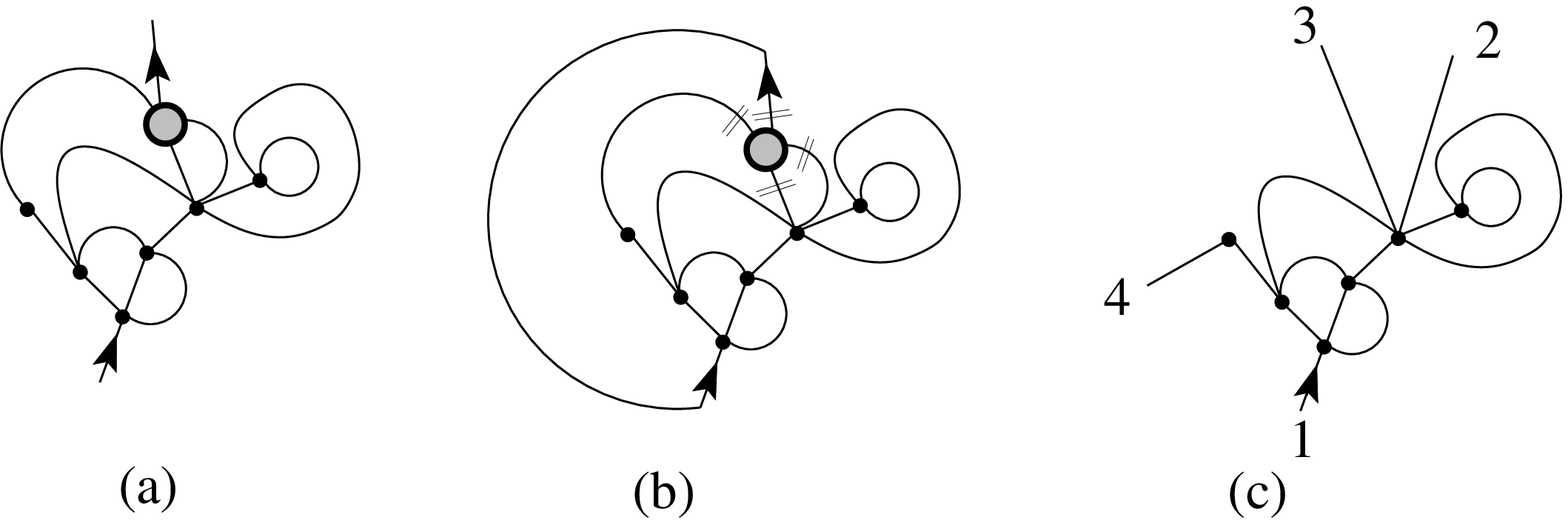}{14.cm}
\figlabel\multipoint
As explained in Sect.2.2, the quantity $Z_{0,-1}(2i-1)$ is the generating function for blossom trees
with a $2i$-valent root vertex and a contour walk with $0$ depth. Returning to the graph 
interpretation, it is identified via the bijection of Sect.2.1 as the generating function
for two-leg diagrams with an outcoming leg attached to a $2i$-valent vertex and 
adjacent to the {\it same} (external) face as the incoming leg. 
The $2i$-valent vertex receives a weight $1$ as opposed
to all other inner vertices weighted by the usual factors $g_k$ if they are $2k$-valent. 
Such two-leg diagrams are in bijection with the $2i$-leg diagrams defined above. Indeed, starting 
from these two-leg diagrams and gluing their two legs into a marked edge, 
we simply erase the unweighted $2i$-valent vertex and obtain the desired $2i$-leg diagrams 
with a marked leg which receives the label $1$ (see Fig.\multipoint\ for an example). 
Note that $2i$-leg diagrams need not be connected except for $i=1$ and may split into connected 
$2j$-leg diagrams with some $j$'s summing to $i$. Eq.\valinteg\ follows.

An outcome of this result is that we obtain a host of explicit expressions for
$G_{2i}$ in terms of the generating functions $R_i$ of blossom trees corresponding
to various choices of $n$ in Eq.\valinteg. The first expression is nothing but
$G_{2i}=Z_{0,-1}(2i-1)$ corresponding to $n=0$. For illustration, for $i=1,2$ and $3$, we have
\eqn\gval{\eqalign{G_{2}& =R_0\cr G_{4}& =R_0(R_0+R_1)\cr G_{6}& 
=R_0 (R_0^2 +2 R_0 R_1 + R_1^2 +R_1 R_2)\cr}}
Another particularly simple choice consists in using Eq.\valinteg\ in the limit $n\to \infty$.
This gives explicit formulas for $G_{2i}$ in terms of the function $R$ solution of Eq.\eqforR\ (recall
that $R$ is simply the generating function for two-leg diagrams with arbitrary geodesic
distance between the legs). Letting $R_a\to R$ for all $a$, we first find that
\eqn\zlim{\eqalign{& Z^{+}_{n-1,n-1+2j}(2i) \to \left({2i\choose i-j}-
{2i\choose i-j-1}\right) R^{i-j}\cr
& V'_{n+2j-1,n-2} \to \sum_{k\geq j+1} g_{k} {2k-1\choose k+j}R^{k+j}\cr}}
which finally yields from Eq.\defgamma
\eqn\resugam{\eqalign{G_{2i}=& \left({2i\choose i}-{2i\choose i-1}\right)R^i\cr 
& -\sum_{k\geq 2} g_k R^{i+k}
\sum_{j=1}^{{\rm min}(i,k-1)} \left({2i\choose i-j}-{2i\choose i-j-1}\right) 
{2k-1\choose k+j}\cr}}
These expressions are determined equivalently as the solutions of the inverse 
relation \Vofgamma\ at large $n$, which reads
\eqn\alterrel{\sum_{k\geq j+1}g_k{2k-1\choose k+j}R^{k+j}=\sum_{i=0}^{j} (-1)^{i-1} 
{2j-i\choose i} R^i G_{2j-2i}}
where ${2j-i \choose i}$ is the number of hard dimer configurations of $i$ dimers 
on a segment of length $2j-1$

\noindent For illustration in the truncated case with valences up to $2m=6$, we find
from Eq.\resugam\ that 
\eqn\resugamsamp{\eqalign{G_2&= R-g_2 R^3 - 5 g_3 R^4\cr G_4&= 2 R^2 -3 g_2 R^4 -16 g_3 R^5
\cr G_6&= 5 R^3-9 g_2 R^5 -50 g_3  R^6 \cr}}
where $R$ is now determined by Eq.\eqforRtrunc.

The explicit expressions \resugam\ or \alterrel\ are to be compared with other
known expressions for multi-point correlation functions. The $2i$-point function above 
may indeed be computed alternatively either by use of the planar limit of the one matrix 
model or by the so-called loop equations. 
The planar solution of the one-matrix integral may be obtained via 
saddle point techniques and reads in the so-called one-cut case:
\eqn\matmod{\eqalign{\sum_{\ell=0}^{i}G_{2i-2\ell}{2\ell\choose \ell}R^\ell& = \lim_{n\to \infty}
\langle n | Q^{2i}(Q-V'(Q)) | n-1 \rangle \cr &= 
{2i+1\choose i}R^i  -\sum_{k\geq 2} g_k R^{i+k}
\sum_{j=1}^{{\rm min}(i,k-1)} {2i+1 \choose i-j} {2k-1\choose k+j} 
\cr}}
This expression is readily equivalent to Eq.\resugam\ by simply noting that
\eqn\checkcomp{\sum_{\ell=0}^{i}{2\ell \choose \ell}\left({2i-2\ell\choose i-\ell-j}-
{2i-2\ell \choose i-\ell-j-1}\right)={2i+1 \choose i-j}}
obtained by cutting any walk of length $2i+1$ from, say $0$ to height $2j+1$ at
the level of its last $0\to 1$ step, resulting into a first walk of length,
say $2\ell$ from $0$ to $0$ and a positive walk of length $2(i-\ell)$ from
height $1$ to height $2j+1$. 

On the other hand, the loop equations simply express the (disconnected) $2i+2$-point function 
as a sum of either a $2i+2k$-point function if the first leg is connected to
a $2k$-valent vertex or a product of two lower order correlations if
the first leg is connected directly to the $(2j+2)$-th leg, namely:
\eqn\loopeq{G_{2i+2}=\sum_{k\geq 1}g_k G_{2i+2k} + \sum_{j=0}^{i} G_{2i-2j} G_{2j}}
That Eq.\resugam\ solves this loop equation may be seen as a necessary consistency of
the planar limit of the one-matrix model in the one-cut case which implies that
Eqs.\loopeq\ and \matmod\ are compatible.

\subsec{Relations between conserved quantities for bounded valences}

In the truncated case just above of valences up to $2m=6$,  Eqs.\resugamsamp\ imply
the relation $g_1 G_2+g_2 G_4 + g_3 G_6 +1 = G_2$ provided Eq.\eqforRtrunc\ is satisfied. 
This may be seen as a particular case of the loop equation \loopeq\ for $i=0$ above in 
its truncated form. More generally, the truncated loop equations for graphs with
valences up to $2m$ allow for expressing all $G_{2i}$ for $i\geq m$ as polynomials
of the first $m-1$ values $G_2,G_4,\ldots,G_{2(m-1)}$. This gives a set of polynomial
relations between the {\it values} of the conserved quantities $\Gamma_{2i}(n)$.

Using the results of previous Section, this interdependence may be rephrased into {\it linear} 
relations between the $G_{2i}$'s by writing Eq.\alterrel\ 
for $j\geq m$ in which case the l.h.s. vanishes. This results in
\eqn\alterrelm{0=\sum_{i=0}^{j} (-1)^{i-1} {2j-i\choose i} R^i G_{2j-2i}}
where $R$ is determined by the polynomial truncation of Eq.\eqforR\  
with $g_k=0$ for $k>m$.
For illustration, for $m=3$, we have the first two relations
\eqn\altertrois{\eqalign{G_6& =  5 R G_4 - 6R^2 G_2 +R^3 \cr 
G_8 &= 7 R G_6 - 15 R^2 G_4 +10 R^3 G_2 - R^4\cr} }
and similar relations for $G_{2i}$'s with higher indices. 

\subsec{Combinatorial interpretation of $\Gamma_{2}(n)$}
\fig{Any blossom tree (a) whose root is encircled by $i\geq 1$ edges in the
closing process of Sect.2.1 may be rerooted at the first excess bud clockwise 
from the root, while this original root is replaced by a leaf (b). This results in a 
rooted tree with a, say $2k$-valent root vertex adjacent to $k-2$ buds and $k+1$ blossom subtrees
(with necessarily $k\geq 2$). The root of this new tree is encircled in
the closing process by $i-1$ edges. In (c) and (d) we focus on the
rearrangement of bud-leaf pairs in a general case. The excess buds in
(d) are those of (c) except the first one, now promoted to root. This
leaves three unmatched leaves.}{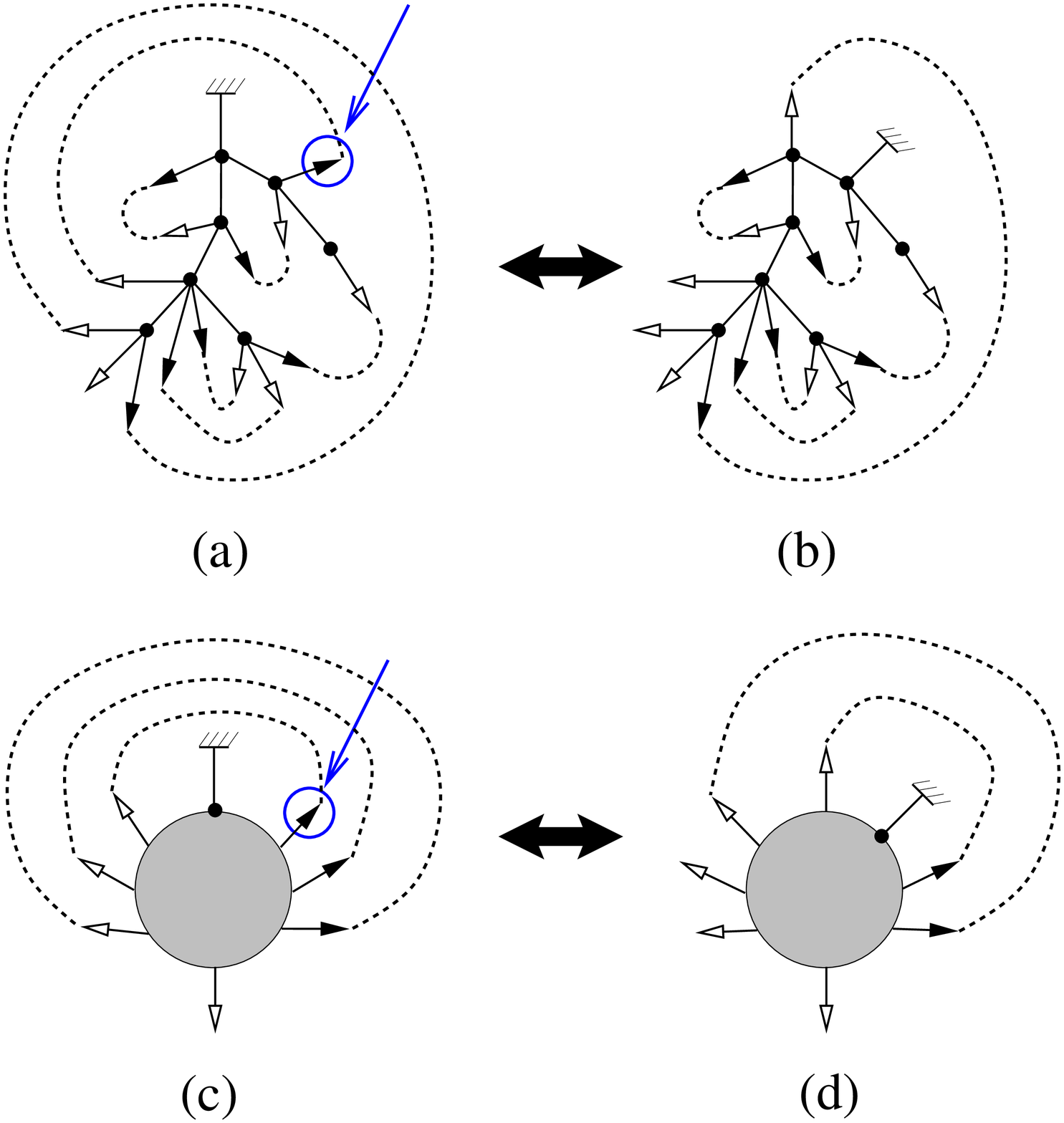}{10.cm}
\figlabel\reroot
The conservation of $\Gamma_2(n)$ has a nice combinatorial explanation in the language of 
blossom trees. Indeed, assuming $\Gamma_0(n)=1$ for all $n\geq 0$, the equality 
$\Gamma_2(n)=\Gamma_2(0)$ for all $n\geq 1$, with $\Gamma_2(n)$ defined by Eq.\defgamma, 
may be rewritten as
\eqn\combidtwo{(R_n-R_0)=V'_{n+1,n-2}}
The quantity $R_n-R_0$ is the generating function for blossom 
trees with contour walk of depth $i$ with $0<i\leq n$. In other words, 
in the closing process of these trees into two-leg diagrams, their root is encircled 
by at least one and at most $n$ edges separating it from the external face. 
Picking the bud from which the deepest encircling edge originates,
we may reroot the tree at this bud and replace the original root by a leaf
(see Fig.\reroot). The resulting object is a tree satisfying (B1), (B2) and (B3) 
except for the root vertex which now has exactly $(k-2)$ buds if it is $2k$-valent 
(note that by construction, this vertex has valence at least $4$ as it originally
carried a bud). In particular, all the proper subtrees of this tree are 
ordinary blossom trees. Finally, the contour walk of this new
tree (now stepping from height $0$ to height $+3$) has clearly depth $(i-1)\leq (n-1)$.
The desired generating function may again be obtained by reading the sequence 
of buds and subtrees counterclockwise around the root vertex and reads
$Z_{n+1,n-2}(2k-1)$ if the new root vertex is $2k$-valent. Summing over
all possible valences, we get the generating function $V'_{n+1,n-2}$. As
the above re-rooting procedure clearly establishes a bijection between the two 
types of trees at hand, this provides a {\it bijective proof} of the equality 
\combidtwo. As this equality is valid for all $n\geq 1$, this in turn proves 
the conservation of $\Gamma_2(n)$.

It would be nice to have similar bijective proofs for the
conservation of $\Gamma_{2i}(n)$ for $i\geq 2$ as well. For $2i=4$, it is
easily seen that, upon using $\Gamma_0(n)=1$ and $\Gamma_2(n+2)=\Gamma_2(n)=R_0$,
the equality $\Gamma_4(n)=\Gamma_4(0)$ may be rewritten as
\eqn\combidfour{R_0(R_{n+1}-R_{1})= V'_{n+3,n-2}+V'_{n+3,n}V'_{n+1,n-2}}
The l.h.s. of this identity is nothing but the generating function for {\it pairs} made of
a blossom tree with contour walk of depth $0$ and a blossom tree with contour walk of 
depth $i$ with $1<i\leq n+1$. The r.h.s. may be interpreted as the generating function for 
the collection of two types of objects
\item{(1)} rooted trees satisfying (B1), (B2) and (B3) except at the root vertex which
carries exactly $(k-3)$ buds if it is $2k$-valent (with necessarily $k\geq 3$). The
contour walk of this tree moreover has depth at most $n-1$;
\item{(2)} pairs of rooted trees satisfying (B1), (B2) and (B3) except at their root vertex which
carries exactly $(k-2)$ buds if it is $2k$-valent (with necessarily $k\geq 2$). The first of
these trees has a contour walk of depth at most $n+1$ and the other a contour walk of
depth at most $n-1$.
\par 
\noindent 
This suggests the existence of a bijection between the two collections of objects above.
If such a bijection could be exhibited, it would provide an alternative 
bijective proof of the conservation of $\Gamma_4(n)$. 

For larger $i$, the conservation of $\Gamma_{2i}(n)$ yields more involved relations
of the type above which may suggest higher order bijections as well. For instance
for $2i=6$, we get the relation
\eqn\combidsix{\eqalign{R_0 R_1(R_{n+2}& -R_{2})- R_0 (R_{n+1}-R_1)(R_{n+3}-R_1)  = 
\cr & V'_{n+5,n-2} +V'_{n+5,n}V'_{n+1,n-2} +
V'_{n+5,n+2}V'_{n+3,n-2}+V'_{n+5,n+2}V'_{n+3,n}V'_{n+1,n-2} \cr}}
while for $2i=8$, we have
\eqn\combideight{\eqalign{& R_0 R_1 R_2(R_{n+3}-R_3)- R_0 R_1\Big( (R_{n+1}-R_1)(R_{n+4}-R_2)
+(R_{n+2}-R_2)(R_{n+4}-R_1) \cr &+(R_{n+2}-R_2)(R_{n+5}-R_2)\Big)+R_0\Big( 
(R_{n+1}-R_1)(R_{n+3}-R_1)(R_{n+5}-R_1)\Big) \cr & = 
V'_{n+7,n-2}+V'_{n+7,n}V'_{n+1,n-2}+V'_{n+7,n+4}V'_{n+5,n-2}\cr & + V'_{n+7,n+2}V'_{n+3,n}V'_{n+1,n-2}
+V'_{n+7,n+4}V'_{n+5,n}V'_{n+1,n-2}+V'_{n+7,n+4}V'_{n+5,n+2}V'_{n+3,n-2}\cr & +V'_{n+7,n+2}V'_{n+3,n-2}
+V'_{n+7,n+4}V'_{n+5,n+2}V'_{n+3,n}V'_{n+1,n-2} \cr}}
These relations still await a bijective proof.

\newsec{Discussion}

\subsec{Alternative expressions for the conserved quantities}
As already mentioned, the definition \defgamma\ corresponds to a 
particular choice of conserved quantities involving positive walks.
This choice however does not respect the ``time reversal" symmetry
of the master equation, namely under $R_{n-j}\leftrightarrow R_{n+j}$
for all $j$ (see Eq.\masterbth\ for illustration). In this Section, 
we present another choice of conserved quantities that respect this
symmetry. More precisely, we define
\eqn\defgamtilde{\encadremath{\eqalign{{\tilde \Gamma}_{2i}(n)&=\left\{ Z_{n,n-1}(2i-1)-R_{n-1} 
R_{n+1} Z_{n-2,n+1}(2i-1)\right\} {\tilde \Gamma}_0(n) \cr &
-\sum_{j=1}^{i} \left\{Z_{n-j,n+j-1}(2i-1)-R_{n-j-1}R_{n+j+1}Z_{n-j-2,n+j+1}(2i-1)\right\} V'_{n+j,n-j-1}\cr}} }
for $i\geq 1$ and $n\geq 0$, with ${\tilde \Gamma}_0(n)$ given by
\eqn\defgamzeroti{{\tilde \Gamma}_0(n)= R_{n}-V'_{n,n-1}}
for $n\geq 0$. Note that in both $Z$'s and $V'$'s, the indices lead to symmetric expressions.

By techniques similar to those of Sect.4.1, one can invert Eqs. \defgamtilde\ and
\defgamzeroti\ into
\eqn\Voftilde{V'_{n+j,n-j-1}=\sum_{i=0}^{j} (-1)^{i-1} \Pi_{n-j,n+j-1}(i) {\tilde \Gamma}_{2j-2i}(n) 
+\delta_{j,0} R_{n}}
with the convention $\Pi_{n,n-1}(0)=1$. The passage from Eq.\defgamtilde\ to Eq.\Voftilde\
makes use of an inversion relation similar to Eq.\geninversion. This relation
is explicited in Appendix B.
 
The conservation of ${\tilde \Gamma}_{2i}(n)$ may be derived along the same
lines as in Sect.4.2.
Using again Eq.\fundV\ now with $a=n+j$ and $b=n-j$, we have the relation
\eqn\fundtilde{\eqalign{V'_{(n+1)+j,(n+1)-j-1}+R_{n+j}& V'_{n+(j-1),n-(j-1)-1} \cr &
=V'_{n+j,n-j-1}+R_{n-j+1}V'_{(n+1)+(j-1),(n+1)-(j-1)-1}\cr}}
For $j\geq 1$ and $n\geq 0$, substituting Eq.\Voftilde\ into this relation yields
\eqn\prooftilde{\eqalign{\big(-{\tilde \Gamma}_{2j}(n+1)& +\sum_{i=1}^{j} (-1)^{i-1} \Pi_{(n+1)-j,(n+1)+j-1}(i) 
{\tilde \Gamma}_{2j-2i}(n+1) \big) \cr 
+R_{n+j} & \big(\sum_{i=0}^{j-1} (-1)^{i-1} \Pi_{n-(j-1),n+(j-1)-1}(i) {\tilde \Gamma}_{2(j-1)-2i}(n) 
+\delta_{(j-1),0} R_{n}\big)\cr =
\big(-{\tilde \Gamma}_{2j}(n)& +\sum_{i=1}^{j} (-1)^{i-1} \Pi_{n-j,n+j-1}(i) {\tilde \Gamma}_{2j-2i}(n) \big) 
\cr
+R_{n-j+1} \big(\sum_{i=0}^{j-1}& (-1)^{i-1} \Pi_{(n+1)-(j-1),(n+1)+(j-1)-1}(i) 
{\tilde \Gamma}_{2(j-1)-2i}(n+1)
+\delta_{(j-1),0} R_{n+1}\big)\cr}}
As the boundary terms on both sides cancel, upon setting $i\to (i-1)$ 
in the last sum of each side, we may rewrite this identity as 
\eqn\reproofti{\eqalign{
{\tilde \Gamma}_{2j}(n+1) & -{\tilde \Gamma}_{2j}(n) \cr
& = \sum_{i=1}^{j} (-1)^{i-1} \left\{ \Pi_{n-j+1,n+j}(i) + R_{n-j+1} 
\Pi_{n-j+2,n+j-1}(i-1) \right\} 
{\tilde \Gamma}_{2j-2i}(n+1) \cr 
& \ \ -\sum_{i=1}^{j} (-1)^{i-1} \left\{ \Pi_{n-j,n+j-1}(i) + R_{n+j} \Pi_{n-j+1,n+j-2}(i-1)\right\} 
{\tilde \Gamma}_{2j-2i}(n) \cr
&= \sum_{i=1}^{j} (-1)^{i-1} \left\{\Sigma_{n-j+1,n+j-1}(i) \right\} \big( {\tilde \Gamma}_{2j-2i}(n+1)- 
{\tilde \Gamma}_{2j-2i}(n)\big)}}
where we have defined 
\eqn\defsig{\Sigma_{a,b}(i)\equiv \Pi_{a,b+1}(i) +\Pi_{a-1,b}(i) -\Pi_{a,b}(i)}
In the last line of Eq.\reproofti\ we have used the relation $\Pi_{a,b+1}(i)+R_{a}\Pi_{a+1,b}(i-1)
= \Pi_{a-1,b}(i)+R_{b+1} \Pi_{a,b-1}(i-1)= \Sigma_{a,b}(i)$ which is a direct consequence of Eq.\eqhardgen. 
We immediately deduce from Eq.\reproofti\ that all ${\tilde \Gamma}_{2j}(n)$ for $j\geq 1$ are conserved 
provided ${\tilde \Gamma}_0(n)$ is independent of $n$. This is the case when the master 
equation \master\ is satisfied, as it reads precisely ${\tilde \Gamma}_0(n)=1$ for all $n\geq 0$.
This completes the proof of conservation of ${\tilde \Gamma}_{2j}$.

For illustration, we list below the two independent conserved quantities for
the truncated case with up to $2m=6$-valent graph corresponding to the 
truncated master equation \masterbth:
\eqn\firsttild{\eqalign{{\tilde \Gamma}_2(n)&= R_n -V'_{n+1,n-2}\cr
{\tilde \Gamma}_4(n)& = R_n(R_{n-1}+R_n+R_{n+1})-R_{n-1}R_{n+1}\cr
 & \ \ \ -(R_{n-1}+R_{n}+R_{n+1})V'_{n+1,n-2}-V'_{n+2,n-3}\cr}}
with
\eqn\firstvpriti{
\eqalign{ V'_{n+1,n-2}& = R_{n+1}R_n R_{n-1}\big(g_2 +g_3 (R_{n+2}+R_{n+1}+R_n+R_{n-1}+R_{n-2})\big)\cr
V'_{n+2,n-3}& = g_3 R_{n+2} R_{n+1} R_{n}R_{n-1} R_{n-2}\cr}}
where we explicitly substituted ${\tilde \Gamma}_0(n)=1$ into Eq.\defgamtilde. 

When compared with Eq.\firstcons, we see that ${\tilde \Gamma}_2(n)=\Gamma_2(n)$
and ${\tilde \Gamma}_4(n)=\Gamma_4(n-1)+(R_{n+1}+R_n+R_{n-1})(\Gamma_2(n)-\Gamma_2(n-1))$.
This shows that the new set of conserved quantities is not independent from that of
Section 4, as expected.

This interdependence can be made even more explicit by noting that, for $n=0$,
we again have ${\tilde \Gamma}_{2i}(0)=Z_{0,-1}(2i-1)=\Gamma_{2i}(0)$. When 
Eq.\master\ is satisfied, the constant value of ${\tilde \Gamma}_{2i}(n)$ is therefore
equal to $G_{2i}$. This is also apparent from the large $n$ limit of Eq.\Voftilde\
which is identical to the large $n$ limit of Eq.\Vofgamma.
As a last remark, we note that the conservation of the $\Gamma$'s is granted by
that of the ${\tilde\Gamma}$'s as we may then replace ${\tilde \Gamma}_{2j-2i}(n)$ in
the r.h.s. of Eq.\Voftilde\ by a shifted value ${\tilde \Gamma}_{2j-2i}(n-j+1)$ which,
upon a global shift of $n\to n+j-1$, show that the ${\tilde\Gamma}_{2i}$'s
also obey Eq.\Vofgamma\ for $j\geq 1$. Together with the master equation which
guarantees that ${\tilde \Gamma}_0(n)=\Gamma_0(n)=1$, this implies that, when
the ${\tilde\Gamma}$'s are conserved, we necessarily have $\Gamma_{2i}(n)={\tilde \Gamma}_{2i}(n)$ 
for all $n$, hence these are conserved as well and take the same constant values.

\subsec{Compacted conserved quantities}

Let us again restrict ourselves to the truncated case of valences up to $2m$ and
concentrate on the $(m-1)$ first conserved quantities ${\tilde \Gamma}_{2i}(n)$
for $i=1,\ldots,m-1$. From the definition \defgamtilde, we note that the largest
index carried by an $R_a$ comes from the contributions $g_m Z_{n+j,n-j-1}(2m-1)$ 
in the term $V'_{n+j,n-j-1}$ appearing in the r.h.s.. This index is 
equal to  $a=n+m-1$ and reached for each $j$ by a unique path of length $2m-1$ starting with
$m-j-1$ up steps followed by $m+j$ down steps. This term reads precisely
$-g_m \prod\limits_{\ell=1-m}^{j}R_{n-\ell}$.
Similarly, the largest index carried by an $R_a$ in ${\tilde \Gamma}_0(n)$
comes from the contribution $-g_m Z_{n,n-1}(2m-1)$ in the term $-V'_{n,n-1}$ 
and reads $-g_m \prod\limits_{\ell=1-m}^{0}R_{n-\ell}$.

We have the remarkable identity 
\eqn\remarkZ{\eqalign{Z_{n-1,n-1}&(2i)=\left\{ Z_{n,n-1}(2i-1)-R_{n-1}
R_{n+1} Z_{n-2,n+1}(2i-1)\right\} \cr &
+\sum_{j=1}^{i} \left\{Z_{n-j,n+j-1}(2i-1)-R_{n-j-1}R_{n+j+1}Z_{n-j-2,n+j+1}(2i-1)\right\} 
\prod_{\ell=1}^j R_{n-\ell}\cr}}
which may be inverted into
\eqn\prodR{-\prod_{\ell=1}^j R_{n-\ell}=\sum_{i=0}^{j} (-1)^{i-1} \Pi_{n-j,n+j-1}(i) 
Z_{n-1,n-1}(2j-2i)}
These identities may be proved with arguments similar to that of Sect.3.2
and follow from general inversion formulas listed in Appendix B.
We deduce that in the combination
\eqn\compac{\encadremath{\Theta_{2i}(n)={\tilde \Gamma}_{2i}(n)+\left(1-{\tilde \Gamma}_0(n)\right) 
Z_{n-1,n-1}(2i)}}
all terms containing $R_{n+m-1}$ are cancelled. Note that the definition \compac\ does
not involve the particular value of $m$ and therefore gives a universal recipe
for compacting the conserved quantities, irrespectively of the precise (finite) order of
truncation.
The quantities $\Theta_{2i}(n)$ form clearly an alternative set of conserved quantities
that take the {\it same} constant values $G_{2i}$, but that involve one less $R_a$ than 
the original $\Gamma_{2i}$ or ${\tilde \Gamma}_{2i}$.

For illustration, for $m=3$, we have 
\eqn\firsttheta{\eqalign{\Theta_2(n)&= R_n+R_{n-1}-R_nR_{n-1}\Big[1-g_1 -g_2 (R_n+R_{n-1}) \cr
&\ \ \ \ \ \ \ \ \ \ \ \ \ \ \ \ \  -g_3 \left(R_{n+1}R_n+(R_n+R_{n-1})^2 
+R_{n-1}R_{n-2}-R_{n+1}R_{n-2}\right)\Big]\cr 
\Theta_4(n)&= R_{n+1}R_n+R_{n-1}R_{n-2}+(R_{n}+R_{n-1})^2 \cr
&\ -R_nR_{n-1}\Big[(1-g_1)(R_{n+1}+R_n+R_{n-1}+R_{n-2}) \cr 
&\ -g_2  (R_{n+1}+R_n+R_{n-1})(R_n+R_{n-1}+R_{n-2})\cr
&\ -g_3 (R_{n+1}+R_n+R_{n-1}+R_{n-2})\left(R_{n+1}R_n+(R_n+R_{n-1})^2+R_{n-1}R_{n-2}\right)\Big]\cr}}
The value of $\Theta_2(n)$ at $g_1=g_3=0$ was already obtained in Ref.\ONEWALL.

\subsec{Application to tetra- and hexavalent graphs}
In the truncated case of graphs with valences up to $2m$, the truncated master equation
is a recursion relation on $n$ allowing for the recursive determination of all
$R_n$'s with $n\geq m-1$ from the knowledge of $R_0,R_1,\ldots,R_{m-2}$. These
initial values are completely determined by using the $m-1$ first conserved quantities
in the form $\Gamma_{2i}(0)=\Gamma_{2i}(\infty)=G_{2i}$ (or similar equalities
for ${\tilde \Gamma}_{2i}$ or $\Theta_{2i}$) with $i=1,\ldots,m-1$.
For instance, in the case $m=3$ of graphs with bi-, tetra- and hexavalent vertices only, 
this yields, by use of the first two lines of \gval\ and \resugamsamp
\eqn\rzerorone{\eqalign{R_0& = R-g_2 R^3 -5 g_3 R^4 \cr 
R_1 &= {R(1 - 6 g_3 R^3 - g_2^2 R^4 - 25 g_3^2 R^6 - g_2
(R^2 + 10 g_3 R^5))\over 1 - g_2 R^2 - 5 g_3 R^3}
\cr}}
with $R$ solution of Eq.\eqforRtrunc. This result is in agreement with
the expressions obtained in Ref.\GEOD\ via the exact solution \sol.

Another possible application concerns local environments within graphs.
For instance, we have access to the distribution of the {\it number of faces
adjacent to the external face} in rooted planar graphs, namely graphs with a 
marked oriented edge, with the external face adjacent to this rooted edge 
and lying on its right. These rooted graphs are in bijection with 
two-leg diagrams with both legs in the same (external) face, as readily seen
by gluing the two legs counterclockwise around the diagram into a rooted edge oriented 
from the outcoming leg to the incoming one. The rooted graphs are therefore 
enumerated by $R_0$. We may then assign an extra weight $x$ to each face adjacent
to the external one, i.e. a weight $x^p$ whenever the external face has $p$ adjacent faces.
The corresponding generating function, still denoted $R_0\equiv R_0(x)$, is
determined by the master equation \master\ for all $n\geq 1$ and
by a modified equation at $n=0$, namely
\eqn\modmast{R_0=x+V'_{0,-1}}
In particular, all conserved quantities are valid for $n\geq 1$ and
their constant value remains unchanged, given by $G_{2i}$.

In the truncated case of valences up to $2m$, we obtain a closed algebraic system
involving $R_0,R_1,\ldots,R_{m-1}$ by supplementing Eq.\modmast\ by the
$m-1$ first conservation laws $\Theta_{2i}(1)=G_{2i}$ for $i=1,\ldots m-1$.
Upon eliminating all $R_n$'s with $n>0$, we are left with a single
algebraic equation determining $R_0(x)$. Let us now illustrate this in
the cases of pure tetravalent and pure hexavalent graphs.
\fig{The twelve rooted tetravalent planar graphs with up to $2$ vertices, hence
with weights $g_2^p$, $p=0,1,2$. Assigning a weight $x$ to each face
adjacent to the external one, we reproduce the first three terms of the
expansion (6.17)}{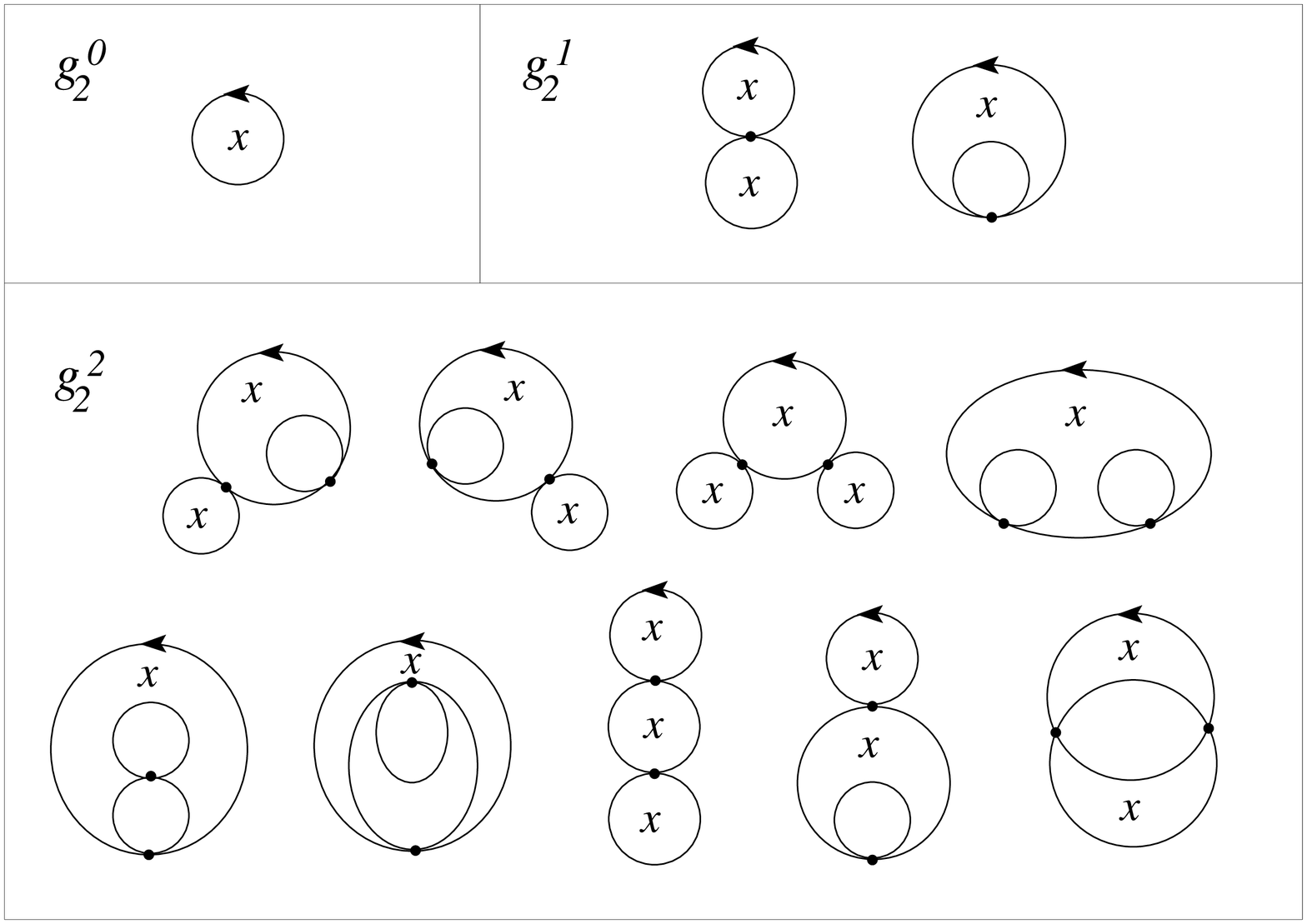}{12.cm}
\figlabel\douze
For pure tetravalent graphs, we obtain the equation
\eqn\Rzerfour{x(x-1) + (1-x + g_2R (R-2) -2 g_2^2 R^3) R_0 + x g_2 R_0^2=0 }
with $R$ solution of $R=1+3 g_2 R^2$.
This leads to the expansion
\eqn\expfour{R_0=x +g_2 (x+x^2) +g_2^2 (3x+4 x^2 +2 x^3) +g_2^3 (14 x+20 x^2 +15 x^3 +5 x^4) + \ldots} 
The corresponding rooted graphs up to order $g_2^2$ are displayed in Fig.\douze\ for illustration.
{}From the knowledge of $R_0(x)$, we may extract the distribution of the number of faces adjacent to
the external face in {\it large tetravalent rooted graphs} by considering the
limit 
\eqn\defproba{\Delta(x)\equiv \sum_{p=1}^\infty x^p {\cal P}(p) = 
\lim_{N\to\infty} {R_0(x)|_{g_2^N}\over R_0(1)|_{g_2^N}}}
Here ${\cal P}(p)$ denotes the probability for the external face of infinitely large 
tetravalent rooted graphs to have $p$ adjacent faces.
The large $N$ limits above are governed by the approach to the critical point
$g_2\to g_2^*= 1/12$ where $R$ and $R_0$ become singular. Expanding in powers 
of $\epsilon=\sqrt{(g_2^*-g_2)/g_2^*}$, the above ratio giving $\Delta(x)$ is
obtained as the ratio of the coefficients of $\epsilon^3$ in the expansions
of $R_0(x)$ and $R_0(1)$ respectively. Using Eq.\Rzerfour, we get the following 
algebraic equation for $\Delta(x)$
\eqn\eqdeltafour{27x (x-1) +4 (4-3x)^3 \Delta (1-x \Delta)=0}
hence
\eqn\expdelfour{\Delta(x)={1\over 2x}\left(1-{8-9x \over (4-3x)^{3\over 2}}\right)}
and
\eqn\resprob{{\cal P}(p)=\left({3\over 16}\right)^{p+1} {(2p+1)!\over (p+1)!(p-1)!}}

For pure hexavalent graphs, we obtain the equation
\eqn\Rzersix{x(x-1)^2-(x-1)(x-1+2 g_3 R^2(3-2R)+24 g_3^2 R^5)R_0+R x g_3(R-2-5 g_3 R^3) R_0^2
+x^2 g_3 R_0^3=0}
with $R$ solution of $R=1+10 g_3 R^3$.
This now leads to the expansion
\eqn\expsix{R_0=x +g_3 (2x +2 x^2 +x^3 ) +g_3^2 (28 x+ 32 x^2+ 25 x^3 + 12 x^4 +3 x^5) + \ldots} 
{}From $R_0(x)$, we may again extract the distribution of the number of faces adjacent to
the external face in {\it large hexavalent rooted graphs} by considering the
same limiting ratio as in \defproba\ with $g_2\to g_3$ and where ${\cal P}(p)$ now
denotes the probability for the external face of infinitely large hexavalent rooted graphs
to have $p$ adjacent faces.
The corresponding large $N$ limits are now governed by the approach to the critical point
$g_3\to g_3^*= 2/135$ where $R$ and $R_0$ become singular. Performing a
similar expansion as before, we arrive at the following 
algebraic equation for $\Delta(x)$
\eqn\eqdeltasix{125 x (x-1)^2 -9 (x-1)(125 x^3 +475 x^2 +200 x -96) \Delta -27 x (x+4) (7-5x)^3 \Delta^2 
(1-x \Delta)=0}
from which we deduce the probabilities 
\eqn\firstprobas{{\cal P}(1)={125\over 864},\quad {\cal P}(2)={1625\over 10368},\quad
{\cal P}(3)={865625\over 5971968}}
for the external face to have $1,2$ or $3$ neighboring faces.

\newsec{Conclusion}

In this paper, we have explicitly constructed various sets of conserved quantities
for the (integrable) master equation which determines the generating function
$R_n$ for two-leg diagrams with vertices of even valences. These conserved quantities
are constructed in terms of generating functions for random walks, whose relations 
with heaps of dimers and hard dimers were instrumental to grant 
the conservation.

The precise properties that we used, namely some boson/fermion correspondences
and their associated inversion relations, as well as some fundamental relations
satisfied by either walks or hard dimers, are very general and presumably
can be adapted to other cases of interest. Indeed, in the planar graph
enumeration framework, it was already recognized that many other classes of
(decorated) graphs have generating functions which satisfy integrable systems
of algebraic master equations. These include graphs with arbitrary even or odd
valences, as well as the so-called $p$-constellations considered in Ref.\CONST. 
In both cases, the two-leg diagram generating functions display a form similar
to that of Eq.\sol\ [\xref\GEOD,\xref\RAMA]. In a particular case of $3$-constellations, 
namely that of graphs with bicolored trivalent vertices (so-called bicubic maps), 
a conserved quantity was already identified in Ref.\ONEWALL. We suspect that 
the above methods can be extended to the even more general case of planar bipartite
graphs, also enumerated (without geodesic distance) by two-matrix models.

Beside the dependence on the geodesic distance, we know from matrix-model
analysis that another integrable structure underlies the topological expansion for
the generating functions of graphs with arbitrary genus. It would be desirable
to relate the two directions of integrability. A general theory of discrete integrable 
equations was built in Ref.\JM\ with (Hirota bilinear) master equations involving more 
than one integer index $n$ (see also \GNR\ for the theory of discrete Painlev\'e
equations). This suggests to look for a larger integrable structure 
that would include both the dependence on (possibly several) geodesic distances as well 
as that on topology.  

As explained in Sect.5.3, the conservation properties have a strong 
``bijective flavor" in the language of blossom trees. A bijective
explanation for the conservation would be very instructive and seems to be within reach
in view if the relations \combidsix\ and \combideight. If such a proof exists,
it must involve only generic properties of blossom trees, or equivalently
of planar graphs. This could allow to trace back the somewhat mysterious generic
appearance of conserved quantities in the context of planar graph enumeration problems to 
these generic properties. 

Finally, we may hope that the conservation properties have a natural interpretation
of their own both in the language of labeled mobiles and in that of spatially extended
branching processes. This is yet to be found. Note in this context the recent appearance
of yet another integrable master equation governing the so-called ``embedded binary
trees" studied in Ref.\EMB.

To conclude, the general case of graphs with vertices of even valence is
known to give access to multicritical transitions points of 2D Quantum Gravity
by proper fine tuning of the parameters $g_i$. As we have now a complete
picture of the conserved quantities, we may use them to describe
these multicritical points. In the scaling limit of large graphs, 
the master equation is known to reduce at such points to an integrable differential equation
related to the KdV hierarchy \GEOD. Our conserved quantities provide explicit 
discrete counterparts of the integrals of motion for these equations. These
integrals are indeed recovered from our expressions in the scaling limit. 
\bigskip
\noindent{\bf Acknowledgments:} 
We acknowledge the support of the EU network on ``Discrete Random Geometry", grant 
HPRN-CT-1999-00161. We thank J\'er\'emie Bouttier for useful discussions and active
interactions at an early stage of this work. We thank Mireille Bousquet-M\'elou 
for enlightening discussions on the combinatorics of heaps of pieces.

\appendix{A}{Relation between heaps and walks}

We consider the generating functions $Z_{a,b}(2k-1)$ for random walks of $2k-1$ steps 
from height $a$ to height $b$ with weight $R_i$ per descent $i\to (i-1)$, and
its truncated version $Z^{+}_{a,b}(2k-1)$ restricted to ``positive" walks with heights 
larger or equal to $a$.
Similarly, we consider the generating function $H_{a,b}(k)$ for pyramids of $k$ dimers
with base $[a,b]$  (with $(b-a)$ odd and positive) with weight $R_i$ per dimer in 
the stripe $i$ except for those in the right projection. We also consider its 
truncated half-pyramid version $H^{+}_{a,b}(k)$ with dimers only in stripes 
$i> a$. In both cases, the truncation simply amounts to taking $R_i\to 0$
for all $i\leq a$.

We wish to prove that
\eqn\toprove{Z_{a,b}(2k-1)=H_{a,b}(k), \qquad Z^{+}_{a,b}(2k-1)=H^{+}_{a,b}(k)}
for all $k\geq 1$ and $(b-a)$ odd and positive.
\fig{Bijection between walks from $a$ to $b>a$ with $2k-1$ steps and pyramids of
base $[a,b]$ with $k$ dimers, as explained in the text. In (a), the times/locations
of pebble droppings are indicated by filled black ellipses. Each pebble stays in place
during some time (indicated by a black horizontal segment) until it is possibly 
picked up. In (b), we put a dimer at each pick-up time/location. This configuration
is completed into a pyramid of base [a,b] in (c). The inverse mapping makes
use of the splitting of dimers into up- or down-pointing triangles (d), whose
sequence along each line determines the steps of the walk.}{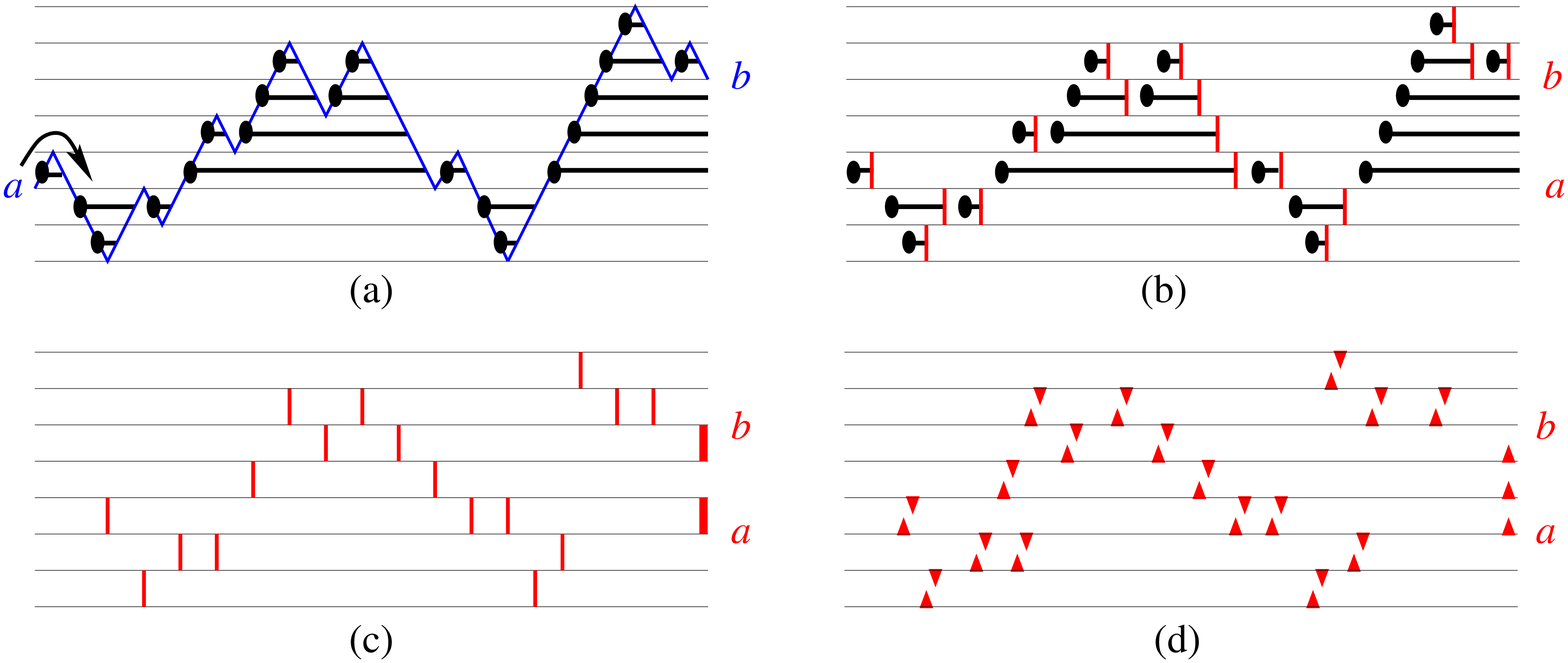}{14.cm}
\figlabel\bijec
These relations are a consequence of the following bijection, borrowed from Ref.\VIENNOT, between on the 
one hand random walks of $2k-1$ steps from height $a$ to height $b>a$ and on the other hand 
pyramids with base $[a,b]$ and a total of $k$ dimers.  
Starting from a walk, as illustrated in Fig.\bijec-(a), viewed as the time (horizontal) evolution of
a walker on the integer (vertical) line, the walker acts as the Little Thumb by
dropping small pebbles at each passage in a unit segment unless a pebble
is already present in which case he picks it up. As shown in Fig.\bijec-(b), putting dimers
at the times/locations of the pick-up events results in a heap of dimers. Moreover, if 
the walker goes from $a$ to $b>a$, there are $b-a$ pebbles left behind. 
By completing the heap on the right by the maximally occupied hard dimer configuration on the
segment $[a,b]$, we obtain a pyramid of base $[a,b]$ and with a total of $k$ dimers 
if the walk has $2k-1$ steps, as illustrated in Fig.\bijec-(c).  
This establishes a bijection whose inverse goes as follows. We split each dimer of
the pyramid not in the right projection into a pair of up- and down-pointing triangles, 
as shown in Fig.\bijec-(d). Each such triangle has its horizontal edge on a line
separating stripes. We also add on the right a succession of $b-a$ up-pointing triangles
adjacent to the lines $y=a,a+1,\ldots,b-1$. For each line $y=i$, we record the
sequence of up- and down-pointing triangles. The walk is reconstructed by starting from
height $a$ and using successively the triangles in the above sequences. More precisely,
at each position $y=i$, the walker makes an up- or down-step according to the 
up-or down-pointing nature of the first unused triangle along the line $y=i$.

In the above bijection, the descending steps of the walk from height $i$ to height
$i-1$ are in one-to-one correspondence with the dimers in the stripe $i$, except
for those of the right projection. Note that these descending steps may correspond 
either to the dropping or to the pick-up of a pebble. In particular, weighting
each descent of the walk from height $i$ to height $i-1$ by $R_i$ (resp. $0$ if 
$i\leq a$) amounts to weight each dimer in the stripe $i$ by $R_i$ (resp. $0$
if $i\leq a$), except for those dimers in the right projection of the pyramid
(resp. half-pyramid). This immediately yields the relations \toprove.

\appendix{B}{Useful identities}
We have the following useful identities, proved by arguments similar to that
of Sect.3.2. 

\noindent For $j\geq 1$ and $k\geq 0$:
\eqn\genuseful{\eqalign{\sum_{i=0}^{j}(-1)^{j-i} \Pi_{n-j,n+j-1}(j-i)
& \Big(\prod_{\ell=1}^k R_{n+k+1-2\ell}\Big) Z_{n-k,n+k-1}(2i-1)\cr
& =\Pi_{n-j,n+j-1}(j) \times \delta_{j\geq k} \times \delta_{j=k\ {\rm mod}\ 2}\cr}}
We immediately deduce from this identity the inversion relation
\eqn\invuse{\eqalign{\sum_{i=0}^{j}(-1)^{j-i}\Pi_{n-j,n+j-1}(j-i)\big(& Z_{n-k,n+k-1}(2i-1)\cr
&\ \ -R_{n-k-1}R_{n+k+1} Z_{n-k-2,n+k+1}(2i-1)\big)=\delta_{j,k}\cr}}
The inversion of Eq.\defgamtilde\ into Eq.\Voftilde\ follows from this identity. 

\noindent For $j\geq 1$:
\eqn\alsouseful{\eqalign{\sum_{i=0}^{j}(-1)^{j-i}\Pi_{n-j,n+j-1}(j-i)
R_{n-1}& Z_{n-2,n-1}(2i-1)\cr &=- \Pi_{n-j,n+j-1}(j) \times \delta_{j=0\ {\rm mod}\ 2} +
\prod_{\ell=1}^j R_{n-\ell}\cr}}
This relation, together with Eq.\genuseful\ for $k=0$, allows to prove Eqs. \remarkZ\ and \prodR\
by noting that $Z_{n-1,n-1}(2i)=Z_{n,n-1}(2i-1)+R_{n-1}Z_{n-2,n-1}(2i-1)$.
\listrefs
\end